\RequirePackage{fix-cm}

\documentclass[smallextended]{svjour3}

\smartqed
\usepackage{graphicx}
\usepackage{mathptmx}
\usepackage{amssymb}

\newcommand{\sotwo}{\mathrm{SO(2)}}
\newcommand{\otwo}{\mathrm{O(2)}}
\newcommand{\sothree}{\mathrm{SO(3)}}
\newcommand{\othree}{\mathrm{O(3)}}

\journalname{Journal of Mathematical Chemistry}

\begin{document}

\title{Molien generating functions and integrity bases for the action of the $\sothree$ and $\othree$ groups on a set of vectors}
\titlerunning{Action of the $\sothree$ and $\othree$ groups on a set of vectors}

\author{Guillaume Dhont \and Patrick Cassam--Chena\"{\i} \and Fr\'ed\'eric Patras}

\institute{Guillaume Dhont \at
  Universit\'e du Littoral C\^ote d'Opale, Laboratoire de Physico--Chimie de l'Atmosph\`ere, MREI2, 189A~avenue Maurice Schumann, 59140 Dunkerque, France \\
  \email{guillaume.dhont@univ-littoral.fr} \\
  ORCID 0000-0001-6184-2971 \\
  \and
  Patrick Cassam--Chena\"{\i} (corresponding author) \and Fr\'ed\'eric Patras \at
  Universit\'e C\^ote d'Azur, CNRS, LJAD, Parc Valrose, 06108 Nice Cedex 2, France \\
  \email{Patrick.Cassam-Chenai@unice.fr, patras@unice.fr}
    ORCID 0000-0002-9437-7794 (pcc) 0000-0002-8098-8279 (fp)
}
  
\date{Received: date / Accepted: date}

\maketitle

\begin{abstract}
The construction of integrity bases for invariant and covariant polynomials built from a set of three dimensional vectors under the $\sothree$ and $\othree$ symmetries is presented. This paper is a follow--up to our previous work that dealt with a set of two dimensional vectors under the action of the $\sotwo$ and $\otwo$ groups [G. Dhont and B. I. Zhilinski\'{\i}, J. Phys. A: Math. Theor., \textbf{46}, 455202 (2013)]. The expressions of the Molien generating functions as one rational function are a useful guide to build integrity bases for the rings of invariants and the free modules of covariants. The structure of the non--free modules of covariants is more complex. In this case, we write the Molien generating function as a sum of rational functions and show that its symbolic interpretation leads to the concept of generalized integrity basis. The integrity bases and generalized integrity bases for $\othree$ are deduced from the $\sothree$ ones. The results are useful in quantum chemistry to describe the potential energy or multipole moment hypersurfaces of molecules. In particular, the generalized integrity bases that are required for the description of the electric and magnetic quadrupole moment hypersurfaces of tetratomic molecules are given for the first time. 
\end{abstract}

\section{\label{S1}Introduction} 

The concept of symmetry pervades many branches of physics and has been formalized by the mathematical framework of group theory \cite{wigner,mcweeny}. In particular, the irreducible representations of a compact Lie group $G$ classify the different ways a set of objects can transform among themselves under the action of the elements $g$ of the group. We define as invariant an object that is unmodified under the elements of the symmetry group $G$, \textit{i.e.} that follows the one dimensional totally symmetric representation $\Gamma_0$ of the group. An irreducible representation $\Gamma \neq \Gamma_0$ of dimension $n$ corresponds to a set of $n$ $\Gamma$--covariants which transform according to the $n \times n$ representation matrix $D^{\left(\Gamma\right)}\left(g\right)$ of the element $g$. Group theory gives methods to build symmetry--adapted objects, however the projector and shift operators \cite{hamermesh} are inefficient tools because the resulting expressions are sometimes redundant or null. Invariant theory \cite{dieudonne,sloane,sturmfels,doc_00628} takes care of the algebraic structure in order to produce powerful methods to construct invariants and covariants. 

Concretely, given a representation $V$ of $G$ (associated for example to $\sothree$ acting on the coordinates of a set of vectors), the action of $G$ extends to the corresponding algebra $A$ of polynomials. An integrity basis (for invariants) is then a set of polynomials that generates the algebra $A_0$ of invariant polynomials. This algebra $A_0$ is Cohen--Macaulay, that is, it is a finitely generated, free module over a polynomial subalgebra: it can be written $A_0=\bigoplus\limits_{i=1}^n{\mathbf C}[g_1,\dots,g_m]\cdot f_i$. The $g_i$ are called primary invariants and form an h.s.o.p. (homogeneous system of parameters), the $f_i$ are called secondary invariants. The decomposition and the h.s.o.p. are not unique but the choice of the h.s.o.p. does not matter for our forthcoming considerations.
To each irreducible representation $\Gamma \neq \Gamma_0$ is associated a finitely generated $A_0$--module of covariants, $A_\Gamma$. If it is free over ${\mathbf C}[g_1,\dots,g_m]$, $A_\Gamma$ is called a Cohen--Macaulay module. When $G$ is finite, all the $A_\Gamma$ are Cohen--Macaulay and the notion of integrity basis for covariants generalizes immediately: it splits into the h.s.o.p. (invariant polynomials) and a family of covariant polynomials that generate freely $A_\Gamma$ as a ${\mathbf C}[g_1,\dots,g_m]$--module. 

The concept of integrity basis has been leveraged in molecular physics \cite{schmelzer,ischtwan,doc_03613,doc_03158,doc_03622,a1776,CCDP_methane}, solid state physics \cite{kopsky1,kopsky2}, physics of deformable bodies \cite{rivlin1,rivlin2,desmorat2020,taurines2020}, high energy physics \cite{sartori}, quantum information \cite{Planat11,Holweck14,Holweck17}, and general multivariate interpolation \cite{RODRIGUEZBAZAN20211}. In recent years, the combination of computer algebraic techniques with artificial intelligence tools to derive primary invariants and fit symmetry--adapted quantum mechanical quantities has been a field of particularly  active research \cite{Li_2013,Shao_2016,Nandi_2019}.

When $G$ is continuous, the Cohen--Macaulay property does not hold in general for $A_\Gamma$ \cite{stanleycov}. The existence of such non--free modules over an h.s.o.p. is one of the noteworthy features unseen when dealing with finite point groups that we want to point out; it requires the introduction of the notion of \textit{generalized integrity basis}. We investigate this phenomenon in detail in the present article.

The present work on the invariants and covariants built from the coordinates of $N$ vectors of the three dimensional space under the $\sothree$ and $\othree$ groups stems from three previous articles \cite{doc_03181,doc_06124,doc_08519}, where the action of the $\sotwo$ and $\otwo$ groups on the coordinates of a set of vectors of the plane was considered. These papers put forward the problem we mentioned of dealing with non--free modules of covariants over an h.s.o.p. that occasionally arises when working with continuous groups. As in our previous paper \cite{doc_03181}, the Molien function \cite{doc_01743} plays a central role in the conception of the integrity bases and their generalization. 

Section~\ref{S2} introduces the setting and presents explicit integrity bases for the rings of invariants and free modules of covariants that are met when dealing with one, two or three vectors. The simplest case of a non--free module occurs for three vectors and irreducible representation $L=2$ and is detailed in Section~\ref{S3}, where we explicitly define and construct the corresponding generalized integrity basis appropriate to the representation of the electric and magnetic quadrupole moment hypersurfaces of tetratomic molecules. We then show that  non--free modules are not uncommon. Section~\ref{S4} presents two conjectures suggested by the results of the two precedent sections, one related to free modules, the other one to non--free modules. They are of practical importance to derive useful bases in the perspective of fitting symmetry--adapted quantum mechanical quantities.

\section{\label{S2}Integrity bases for rings of invariants and free modules of covariants} 

\subsection{\label{S2A}Context and Molien generating functions}

In a polyatomic molecule with $N+1$ nuclei, it is always possible to define $N$ relative vectors that are invariant under any translation \cite{sutcliffe}. These relative vectors can be Jacobi vectors, Radau vectors \cite{GATTI20091} or simply differences of the vector positions of two nuclei, $\vec{r}_i - \vec{r}_j$. From now on, by \textit{vector} we mean one of these relative vectors. 

Each value of the angle of rotation $\omega \in \left[0,\pi\right]$ defines an equivalence class of the $\sothree$ group, with all the rotations of the same angle but different rotation axes $\hat{n}$ belonging to the same class. The character of the class for the $\left(2L+1\right)$--dimensional irreducible representation $\left(L\right)$, $L\in \mathbb{N}$, of $\sothree$ is given by \cite{VMK}:
$$
\chi^{\left(L\right)}\left(\omega\right)
=
\frac{\sin \left[\left(2L+1\right)\frac{\omega}{2}\right]}{\sin \frac{\omega}{2}}.
$$
The three coordinates of one vector span the irreducible representation $\left(1\right)$ of $\sothree$ and the $3N$ coordinates of the $N$ vectors span the reducible representation $\underbrace{\left(1\right) \oplus\cdots \oplus \left(1\right)}_{N \mathrm{\ times}}$, called from now on the initial representation. The group $\sothree$ acts diagonally on the direct sum of the various symmetric tensor powers of this representation or, equivalently, on polynomials in the coordinates. The Molien generating function associated to a given irreducible representation $(L)$ (called the final representation) is then the formal power series in a variable $\lambda$ where the coefficient of $\lambda^n$ is the multiplicity of $(L)$ in the $\sothree$--module of degree $n$ polynomials in the coordinates.

The Molien generating function for computing the number of linearly independent invariants (\textit{i.e.} "$\left(0\right)$--covariants") or $\left(L\right)$--covariants ($L>0$) of a given degree built from $N$ vectors is given by
\begin{equation}
\label{molien1}
g_\sothree\left(N,L;\lambda\right)
=
\frac{1}{2\pi^2}
\int_{0}^{\pi} \int_{0}^{2\pi} \int_{0}^{\pi}
\frac{\chi^{\left(L\right)}\left(\omega\right)}{\det \left(1_{3N}-\lambda
  D\left(\omega,\hat{n}\right)\right)}
\sin \theta d\theta d\varphi \sin^2 \frac{\omega}{2}d\omega,
\end{equation}
where $1_{3N}$ is the $3N \times 3N$ identity matrix and $D\left(\omega,\hat{n}\right)$ is the $3N \times 3N$ representation matrix of the rotation. The angles $\theta$ and $\varphi$ give the orientation of the rotation axis $\hat{n}$. The evaluation of the integral~(\ref{molien1}) is derived in Appendix~\ref{appA}. Collins and Parsons determined the Molien generating function for the invariants with a parametrization of the rotations based on Euler angles \cite{a1610}. Our approach with the rotation angle $\omega$ and rotation axis $\hat{n}$, while giving the same result for the invariants, is more amenable to a generalization for covariants because our parametrization gives a straightforward integral over the single variable $\omega$ while calculations would not be so direct using the Euler angles.

It is remarkable that the Molien generating function for $\sothree$ can be written as the difference between Molien generating functions for $\sotwo$ that were determined in Ref. \cite{doc_03181}:
\begin{eqnarray*}
g_\sothree\left( N,L;\lambda \right)
&=&
\frac{1}{(1-\lambda)^N}
\left(
g_\sotwo\left( N,L ; \lambda \right)
-
g_\sotwo\left( N,L+1 ; \lambda \right)
\right),
\\
&=&
\frac{1}{(1-\lambda)^N}
\left\{
\frac{1}{\pi}
\int_{0}^{\pi}
\frac{\cos \left(L \omega\right)\, d\omega}{(1-2\lambda\cos \omega+\lambda^2)^N}
-\frac{1}{\pi}
\int_{0}^{\pi}
\frac{\cos \left[\left(L+1\right) \omega\right]\, d\omega}{(1-2\lambda\cos \omega+\lambda^2)^N}
\right\},
\nonumber
\end{eqnarray*}
where $g_\sotwo\left( N,L ; \lambda \right)$ is the generating function for the number of covariants with irreducible representation $\left( L \right)$ of $\sotwo$ that can be built from $N$ two dimensional vectors.

\subsection{\label{S2B}Extension to the $\othree$ group}

The matrix representation of the inversion operation for the problem of $N$ vectors is $-1_{3N}$. Denoting the irreducible representations of $\othree$ by $\left(L^\epsilon\right)$, $L\in \mathbb{N}$, $\epsilon=\pm 1$, the corresponding Molien functions for $N$ vectors under the $\othree$ group are obtained through:
$$
g_{\othree}\left( N,L^\epsilon ;\lambda \right)
= \frac{1}{2} \left[ g_\sothree\left( N,L ;\lambda \right)
  +\epsilon g_\sothree\left( N,L ;-\lambda\right) \right].
$$

\subsection{\label{S2C}Integrity bases for the ring of invariants and free modules of covariants}

\subsubsection{\label{S2C1}One vector}

The Molien function for the initial representation $\left( 1 \right)$ resulting from the three coordinates $\mathbf{x}_1=(x_1,y_1,z_1)$ of one vector is
\begin{equation} \label{gen_fun_for_one_vector}
g_\sothree\left(1,L;\lambda \right) = \frac{\lambda^L}{1-\lambda^2}.
\end{equation}
The formal expansion of this generating function gives $\lambda^L+\lambda^{L+2}+\cdots$, which means that we can built one set of $2L+1$ $\left(L\right)$--covariants of degree $p \geq L$ when $p-L$ is even. Formula~(\ref{gen_fun_for_one_vector}) indicates that the polynomials belonging to the free module of $(L)$--covariants can be univoquely built from one primary invariant of degree $2$ and one set of $\left(2L+1\right)$ secondary $\left(L\right)$--covariants of degree $L$. The primary invariant is of course the scalar product $Q_{1,1}=x_1^2+y_1^2+z_1^2$. The $\left(2L+1\right)$ secondary $\left(L\right)$--covariants are naturally chosen as the real solid harmonics $\tilde{\mathcal{Y}}_{L,M}\left(x_1,y_1,z_1\right)$. Appendix~\ref{realSH} gives the expressions of the real solid harmonics for $L$ up to $3$.

\subsubsection{\label{S2C2}Two vectors}

The six coordinates $(\mathbf{x}_1,\mathbf{x}_2)$ of two vectors span the reducible six--dimensional representation $\left(1\right) \oplus \left(1\right)$. The coupling according to formula~(\ref{formule_article_a0628}) of the generating functions for one vector, given in~(\ref{gen_fun_for_one_vector}),  produces the Molien function for two vectors in the two variables $\lambda_1$ and $\lambda_2$:
$$
g_\sothree\left(2,L; \lambda_1, \lambda_2\right)
=
\frac{\sum\limits_{i=0}^L\lambda_1^i\lambda_2^{L-i}+\lambda_1\lambda_2\sum\limits_{i=0}^{L-1}\lambda_1^i\lambda_2^{L-i-1}}{\left(1-\lambda_1^2\right)\left(1-\lambda_2^2\right)\left(1-\lambda_1\lambda_2\right)},
$$
with the convention that the second term in the numerator is zero for $L=0$. If we do not distinguish between the two vectors, we can set $\lambda_1= \lambda_2=\lambda$ to get the Molien function as
\begin{equation} \label{ML2}
g_\sothree\left( 2,L ; \lambda \right)
=
\frac{(L+1)\lambda^L+L\lambda^{L+1}}{\left(1-\lambda^2\right)^3}.
\end{equation}
The coefficients $L+1$ and $L$ in the numerators are always non--negative, suggesting that the expression~(\ref{ML2}) of the Molien function can be used to construct the integrity basis for the invariants and for any $\left(L\right)$--covariants.

The Molien generating function for the ring of invariants ($L=0$) specializes to
$$
g_\sothree\left(2,0;\lambda_1,\lambda_2 \right)
=
\frac{1}{\left(1-\lambda_1^2\right)\left(1-\lambda_2^2\right)\left(1-\lambda_1\lambda_2\right)}.
$$
The ring of invariant is well--known \cite{b100} and the integrity basis consists in the $3$ scalar products $Q_{1,1}=x_1^2+y_1^2+z_1^2$, $Q_{2,2}=x_2^2+y_2^2+z_2^2$, $Q_{1,2}=x_1x_2+y_1y_2+z_1z_2$ as the primary polynomials of degree $2$ and $1$ as the secondary polynomial of degree $0$. The totally invariant quantities such as potential energy hypersurfaces of triatomic molecules can be expanded as polynomials in these three primary polynomials.

For $L=1$, the Molien generating function in two variables is
$$
g_\sothree\left(2,1;\lambda_1,\lambda_2 \right)
=
\frac{\lambda_1+\lambda_2+\lambda_1\lambda_2}{\left(1-\lambda_1^2\right)\left(1-\lambda_2^2\right)\left(1-\lambda_1\lambda_2\right)}.
$$
The three primary polynomials are chosen to be identical to the three scalar products $Q_{1,1}$, $Q_{2,2}$ and $Q_{1,2}$ chosen for the ring of invariants. The two sets of $\left(1\right)$--covariants of degree $1$ are taken as $\sqrt{\frac{4\pi}{3}}\tilde{\mathcal{Y}}_{1,M}\left(\mathbf{x}_i\right)$ with $i \in \left\{1,2\right\}$ and $M \in \left\{1,0,-1\right\}$. The set of $\left(1\right)$--covariants of degree $2$ can be chosen as their cross--product:
$$
\left(\begin{array}{c}
f_{1,1}\left(\mathbf{x}_1,\mathbf{x}_2\right) \\
f_{1,0}\left(\mathbf{x}_1,\mathbf{x}_2\right) \\
f_{1,-1}\left(\mathbf{x}_1,\mathbf{x}_2\right)
\end{array}\right)
= \left( \begin{array}{c}
y_1z_2-z_1y_2\\ x_1y_2-y_1x_2 \\ z_1x_2-x_1z_2
\end{array} \right).
$$
(The subscripts of the $f$--symbols are respectively $L$ and $M$.)

The electric dipole moment function, $\vec{\mu}^\mathrm{el}$, transforms as the irreducible representation $\left(1^-\right)$ of the group $\othree$. The Molien generating function is:
$$
g_{\othree}\left(2,1^-;\lambda_1,\lambda_2\right) =
\frac{\lambda_1+\lambda_2}{\left(1-\lambda_1^2\right)\left(1-\lambda_2^2\right)\left(1-\lambda_1\lambda_2\right)}.
$$
The cross--products do not enter the integrity basis because they are invariant with respect to the inversion operation. Therefore, the components $\mu_M^\mathrm{el}$, $M \in \left\{1,0,-1\right\}$ of the electric dipole moment function of an $\mathrm{ABC}$ molecule have the form:
\begin{equation}
\mu_M^\mathrm{el}\left(\mathbf{x}_1,\mathbf{x}_2\right) = P_1^\mathrm{el}\left(Q_{1,1},Q_{2,2},Q_{1,2}\right) \times \sqrt{\frac{4\pi}{3}}\tilde{\mathcal{Y}}_{1,M}\left(\mathbf{x}_1\right)
+ P_2^\mathrm{el}\left(Q_{1,1},Q_{2,2},Q_{1,2}\right) \times \sqrt{\frac{4\pi}{3}}\tilde{\mathcal{Y}}_{1,M}\left(\mathbf{x}_2\right),
\label{EDMS}
\end{equation}
where $P_1^\mathrm{el}$ and $P_2^\mathrm{el}$ are two polynomials in the three primary invariants $Q_{1,1}$, $Q_{2,2}$, and $Q_{1,2}$. Thus, only two polynomials in three variables need to be fitted with respect to quantum chemistry data to determine the dipole moment surface functions. In contrast, the magnetic dipole moment function,  $\vec{\mu}^\mathrm{mag}$, transforms as the irreducible representation $\left(1^+\right)$ of the group $\othree$, so its components have the form
$$
\mu_M^\mathrm{mag}\left(\mathbf{x}_1,\mathbf{x}_2\right) = P_1^\mathrm{mag}\left(Q_{1,1},Q_{2,2},Q_{1,2}\right) \times f_{1,M}\left(\mathbf{x}_1,\mathbf{x}_2\right),
$$
where $P_1^\mathrm{mag}$ is a polynomial in the three primary invariants.

For $L=2$, the Molien generating function in the two variables is:
$$
g_\sothree^{\left(1\right)}\left(2,2;\lambda_1,\lambda_2\right) =
\frac{\lambda_1^2+\lambda_1\lambda_2+\lambda_2^2+\lambda_1^2\lambda_2+\lambda_1\lambda_2^2}{\left(1-\lambda_1^2\right)\left(1-\lambda_2^2\right)\left(1-\lambda_1\lambda_2\right)}.
$$
The two sets of $(2)$--covariants of degree $2$ corresponding to the polynomials $\lambda_1^2$ and $\lambda_2^2$ in the numerator are chosen as $\sqrt{\frac{4\pi}{5}}\tilde{\mathcal{Y}}_{2,M}\left(\mathbf{x}_i\right)$, $i \in \left\{1,2\right\}$, $M \in \left\{2,1,0,-1,-2\right\}$. The last set of $(2)$--covariants of degree $2$, corresponding to the polynomial $\lambda_1 \lambda_2$, is chosen as:
$$
\left(\begin{array}{c}
f_{2,2}^{(2)}\left(\mathbf{x}_1,\mathbf{x}_2\right) \\
f_{2,1}^{(2)}\left(\mathbf{x}_1,\mathbf{x}_2\right) \\
f_{2,0}^{(2)}\left(\mathbf{x}_1,\mathbf{x}_2\right) \\
f_{2,-1}^{(2)}\left(\mathbf{x}_1,\mathbf{x}_2\right) \\
f_{2,-2}^{(2)}\left(\mathbf{x}_1,\mathbf{x}_2\right)
\end{array}\right)
:=
\left( \begin{array}{c}
\frac{\sqrt{3}}{2}\left(x_1x_2-y_1y_2\right) \\
\frac{\sqrt{3}}{2}\left(x_1z_2+z_1x_2\right) \\
\frac{1}{2}\left(2z_1z_2-x_1x_2-y_1y_2\right) \\
\frac{\sqrt{3}}{2}\left(y_1z_2+z_1y_2\right) \\
\frac{\sqrt{3}}{2}\left(x_1y_2+y_1x_2\right)
\end{array} \right).
$$
The two sets of $(2)$--covariants of degree three, corresponding to the  polynomials $\lambda_1^2\lambda_2$ and $\lambda_1\lambda_2^2$ can be constructed by substituting one set of coordinates in the last expression by the $(L=1)$--cross--product covariant,
$$
\left(\begin{array}{c}
f_{2,2}^{(3,1)}\left(\mathbf{x}_1,\mathbf{x}_2\right) \\
f_{2,1}^{(3,1)}\left(\mathbf{x}_1,\mathbf{x}_2\right) \\
f_{2,0}^{(3,1)}\left(\mathbf{x}_1,\mathbf{x}_2\right) \\
f_{2,-1}^{(3,1)}\left(\mathbf{x}_1,\mathbf{x}_2\right) \\
f_{2,-2}^{(3,1)}\left(\mathbf{x}_1,\mathbf{x}_2\right)
\end{array}\right)
:=
\left( \begin{array}{c}
\frac{\sqrt{3}}{2}\left(x_1(y_1z_2-z_1y_2)-y_1(z_1x_2-x_1z_2)\right) \\
\frac{\sqrt{3}}{2}\left(x_1(x_1y_2-y_1x_2)+z_1(y_1z_2-z_1y_2)\right) \\
\frac{1}{2}\left(2z_1(x_1y_2-y_1x_2)-x_1(y_1z_2-z_1y_2)-y_1(z_1x_2-x_1z_2)\right) \\
\frac{\sqrt{3}}{2}\left(y_1(x_1y_2-y_1x_2)+z_1(z_1x_2-x_1z_2)\right) \\
\frac{\sqrt{3}}{2}\left(x_1(z_1x_2-x_1z_2)+y_1(y_1z_2-z_1y_2)\right)
\end{array} \right),
$$
and
$$
\left(\begin{array}{c}
f_{2,2}^{(3,2)}\left(\mathbf{x}_1,\mathbf{x}_2\right) \\
f_{2,1}^{(3,2)}\left(\mathbf{x}_1,\mathbf{x}_2\right) \\
f_{2,0}^{(3,2)}\left(\mathbf{x}_1,\mathbf{x}_2\right) \\
f_{2,-1}^{(3,2)}\left(\mathbf{x}_1,\mathbf{x}_2\right) \\
f_{2,-2}^{(3,2)}\left(\mathbf{x}_1,\mathbf{x}_2\right)
\end{array}\right)
:=
\left( \begin{array}{c}
\frac{\sqrt{3}}{2}\left((y_1z_2-z_1y_2)x_2-(z_1x_2-x_1z_2)y_2\right) \\
\frac{\sqrt{3}}{2}\left((y_1z_2-z_1y_2)z_2+(x_1y_2-y_1x_2)x_2\right) \\
\frac{1}{2}\left(2(x_1y_2-y_1x_2)z_2-(y_1z_2-z_1y_2)x_2-(z_1x_2-x_1z_2)y_2\right) \\
\frac{\sqrt{3}}{2}\left((z_1x_2-x_1z_2)z_2+(x_1y_2-y_1x_2)y_2\right) \\
\frac{\sqrt{3}}{2}\left((y_1z_2-z_1y_2)y_2+(z_1x_2-x_1z_2)x_2\right) \\
\end{array} \right).
$$
(The superscripts of the $f$--symbols are respectively the degree of the covariant  and an integer index when there are several secondary covariants of the same degree.)

The electric quadrupole moment function,  $q^\mathrm{el}$, transforms as the irreducible representation $\left(2^+\right)$ of the group $\othree$. Its components have the form
\begin{eqnarray}
q_M^\mathrm{el}\left(\mathbf{x}_1,\mathbf{x}_2\right)&=& P_1^\mathrm{el}\left(Q_{1,1},Q_{2,2},Q_{1,2}\right) \times \sqrt{\frac{4\pi}{5}}\tilde{\mathcal{Y}}_{2,M}\left(\mathbf{x}_1\right)+ P_2^\mathrm{el}\left(Q_{1,1},Q_{2,2},Q_{1,2}\right) \times \sqrt{\frac{4\pi}{5}}\tilde{\mathcal{Y}}_{2,M}\left(\mathbf{x}_2\right)
\nonumber \\
&&+P_3^\mathrm{el}\left(Q_{1,1},Q_{2,2},Q_{1,2}\right)  \times f_{2,M}^{(2)}\left(\mathbf{x}_1,\mathbf{x}_2\right),
\label{EQMS}
\end{eqnarray}
while its magnetic counterpart,  $q^\mathrm{mag}$, transforms as the irreducible representation $\left(2^-\right)$ and its components have the form 
\begin{equation}
q_M^\mathrm{mag}\left(\mathbf{x}_1,\mathbf{x}_2\right) = P_1^\mathrm{mag}\left(Q_{1,1},Q_{2,2},Q_{1,2}\right)  \times f_{2,M}^{(3,1)}\left(\mathbf{x}_1,\mathbf{x}_2\right)+P_2^\mathrm{mag}\left(Q_{1,1},Q_{2,2},Q_{1,2}\right)  \times f_{2,M}^{(3,2)}\left(\mathbf{x}_1,\mathbf{x}_2\right).
\label{MQMS}
\end{equation}

\subsubsection{\label{S2C3}Three vectors}

The nine components of the three vectors span the reducible initial representation $\left(1\right) \oplus \left(1\right) \oplus \left(1\right)$ and the Molien generating function in the form of a single rational function with one variable $\lambda$ is:
\begin{equation} \label{g_three_vectors_1}
g_\sothree^{\left(a\right)}\left(3,L;\lambda \right) =
\frac{\mathcal{N}_\sothree^{\left(a\right)}\left(3,L;\lambda\right)}{(1-\lambda^2)^6}, 
\end{equation}
with the numerator
\begin{equation} \label{g_three_vectors_1_}
\mathcal{N}_\sothree^{\left(a\right)}\left( 3,L ; \lambda\right) =
\frac{(L+2)(L+1)}{2}\lambda^L +(L+2)L\lambda^{L+1}
-(L+1)(L-1)\lambda^{L+3} -\frac{L(L-1)}{2}\lambda^{L+4}.
\end{equation}
The coefficients in the numerator~(\ref{g_three_vectors_1_}) are all non--negative if $L=0$ or $L=1$, while negative coefficients appear for $L\geq 2$. This means that the symbolic interpretation of the numerator coefficients as the amount in the integrity basis of $\left(L\right)$--covariants of each degree does not hold anymore. Instead, the first negative coefficient indicates the number of syzygies of lowest degree and diagnoses the fact that the $(L)$--component of the ring of polynomials over the ring of  primary invariants is not a free module. We shall expect to be able to construct an integrity basis for $L=0$ and $L=1$, but the cases with $L \geq 2$ require more consideration.

The Molien generating function for the invariant polynomials ($L=0$) is
$$
g_\sothree^{\left(a\right)}\left(3,0;\lambda_1,\lambda_2,\lambda_3\right) =
\frac{1+\lambda_1\lambda_2\lambda_3}{\left(1-\lambda_1^2\right)\left(1-\lambda_2^2\right)\left(1-\lambda_3^2\right)\left(1-\lambda_1\lambda_2\right)\left(1-\lambda_1\lambda_3\right)\left(1-\lambda_2\lambda_3\right)}.
$$
The $6$ primary invariant polynomials of degree $2$ are chosen as the $6$ scalar products $Q_{i,j}:=x_ix_j+y_iy_j+z_iz_j$ for $1\leq i\leq j\leq 3$ \cite{b100}. The secondary invariant of degree zero is the polynomial constant $1$ and the secondary invariant of degree three is the determinant $\left|\begin{array}{ccc} x_1 & x_2 & x_3\\ y_1 & y_2 & y_3 \\ z_1 & z_2 & z_3 \end{array}\right|$. The determinant transforms as the irreducible representation $\left(0^-\right)$ of the group $\othree$. The potential energy hypersurfaces of tetratomic molecules being totally symmetrical under the  $\othree$ group, they can be expanded as polynomials in the primary invariant polynomials only.

The Molien generating function for $L=1$ is
$$
g_\sothree^{\left(a\right)}\left(3,1;\lambda_1,\lambda_2,\lambda_3\right) =
\frac{\lambda_1+\lambda_2+\lambda_3+\lambda_1\lambda_2+\lambda_1\lambda_3+\lambda_2\lambda_3}{\left(1-\lambda_1^2\right)\left(1-\lambda_2^2\right)\left(1-\lambda_3^2\right)\left(1-\lambda_1\lambda_2\right)\left(1-\lambda_1\lambda_3\right)\left(1-\lambda_2\lambda_3\right)}.
$$
The three sets of $(1)$--covariants of degree $1$ are taken as $\sqrt{\frac{4\pi}{3}}\tilde{\mathcal{Y}}_{1,M}\left(\mathbf{x}_i\right)$, $i\in \left\{1,2,3\right\}$. The three sets of $\left(1\right)$--covariants of degree $2$ are the cross--products $f_{1,M}\left(\mathbf{x}_1,\mathbf{x}_2\right)$, $f_{1,M}\left(\mathbf{x}_1,\mathbf{x}_3\right)$, and $f_{1,M}\left(\mathbf{x}_2,\mathbf{x}_3\right)$. The electric dipole moment function of a tetratomic molecule should thus be expanded according to the form:
\begin{eqnarray}
\mu_M^\mathrm{el}\left(\mathbf{x}_1,\mathbf{x}_2,\mathbf{x}_3\right)
&=&
P_1^\mathrm{el}\left(Q_{1,1},Q_{1,2},Q_{1,3},Q_{2,2},Q_{2,3},Q_{3,3}\right) \times \sqrt{\frac{4\pi}{3}}\tilde{\mathcal{Y}}_{1,M}\left(\mathbf{x}_1\right)
\nonumber \\
&&+ P_2^\mathrm{el}\left(Q_{1,1},Q_{1,2},Q_{1,3},Q_{2,2},Q_{2,3},Q_{3,3}\right) \times \sqrt{\frac{4\pi}{3}}\tilde{\mathcal{Y}}_{1,M}\left(\mathbf{x}_2\right)
\nonumber \\
&&+ P_3^\mathrm{el}\left(Q_{1,1},Q_{1,2},Q_{1,3},Q_{2,2},Q_{2,3},Q_{3,3}\right) \times \sqrt{\frac{4\pi}{3}}\tilde{\mathcal{Y}}_{1,M}\left(\mathbf{x}_3\right),
\nonumber
\end{eqnarray}
where $P_1^\mathrm{el}$, $P_2^\mathrm{el}$, and $P_3^\mathrm{el}$ are polynomials in the $6$ primary invariants $Q_{i,j}$. Similarly, the magnetic dipole moment function is written as: 
\begin{eqnarray}
\mu_M^\mathrm{mag}\left(\mathbf{x}_1,\mathbf{x}_2,\mathbf{x}_3\right)
&=& P_1^\mathrm{mag}\left(Q_{1,1},Q_{1,2},Q_{1,3},Q_{2,2},Q_{2,3},Q_{3,3}\right) \times f_{1,M}\left(\mathbf{x}_2,\mathbf{x}_3\right)
\nonumber \\
&&+ P_2^\mathrm{mag}\left(Q_{1,1},Q_{1,2},Q_{1,3},Q_{2,2},Q_{2,3},Q_{3,3}\right) \times f_{1,M}\left(\mathbf{x}_1,\mathbf{x}_3\right)
\nonumber \\
&&+ P_3^\mathrm{mag}\left(Q_{1,1},Q_{1,2},Q_{1,3},Q_{2,2},Q_{2,3},Q_{3,3}\right) \times f_{1,M}\left(\mathbf{x}_1,\mathbf{x}_2\right).
\nonumber
\end{eqnarray}

\section{\label{S3}Generalized integrity bases for non--free modules}

\subsection{\label{S3A}The simplest case: the $N=3$, $L=2$ non--free module}

We present the need to introduce generalized integrity bases with the simplest possible example: the non--free module for three vectors and $L=2$. When the three representations associated with the three vectors are distinguished by the three variables $\lambda_i$, the numerator of the Molien function is:
\begin{eqnarray}
&\mathcal{N}_\sothree^{\left(a\right)}\left(3,2;\lambda_1,\lambda_2,\lambda_3\right)
=
\lambda_1^2+\lambda_2^2+\lambda_3^2+\lambda_1\lambda_2+\lambda_1\lambda_3+\lambda_2\lambda_3+&\nonumber\\
&\lambda_1^2\lambda_2+\lambda_1^2\lambda_3+\lambda_1\lambda_2^2+\lambda_1\lambda_3^2+\lambda_2^2\lambda_3+\lambda_2\lambda_3^2+2\lambda_1\lambda_2\lambda_3&\nonumber\\
&-\lambda_1^2\lambda_2^2\lambda_3-\lambda_1^2\lambda_2\lambda_3^2-\lambda_1\lambda_2^2\lambda_3^2-\lambda_1^2\lambda_2^2\lambda_3^2.
\label{N_three_vectors_1_L2} 
\end{eqnarray}
As discussed in Section~\ref{S2C3}, the negative coefficients in this numerator make it not suitable for an interpretation of the generating function in term of an integrity basis. However, its symbolic interpretation in term of generators and syzygies suggests the existence of $6$ sets of $\left(2\right)$--covariant generators of degree $2$ and $8$ sets of $\left(2\right)$--covariant generators of degree $3$. These generators are involved in three syzygies of degree $5$ and one syzygy of degree $6$.  

The Molien generating function (\ref{g_three_vectors_1}) for three vectors does not admit a symbolic interpretation in term of an integrity basis for $L \geq 2$ when written as a single rational function. It can nevertheless be recast as a sum of two rational functions where all the coefficients in the numerators are now positive coefficients for $L\geq 2$:
\begin{equation} \label{g_three_vectors_2}
g_\sothree^{\left(b\right)}\left(3,L;\lambda \right) =
\frac{(2L+1)\lambda^L + (2L+1)\lambda^{L+1}}{(1-\lambda^2)^6}
 +\frac{\frac{L(L-1)}{2}\lambda^L + (L+1)(L-1)\lambda^{L+1} + 
\frac{L(L-1)}{2}\lambda^{L+2}}{(1-\lambda^2)^5}.
\end{equation}
The denominators in the two rational functions are different. While the leftmost rational function suggests $6$ primary invariants of degree $2$ as in the $L=0$ or $L=1$ cases, the second rational function points to only five primary invariants of degree $2$ to be selected among the $6$ chosen for the first rational function. The symbolic interpretation of the first rational function specifies that only $2L+1$ sets of secondary $\left(L\right)$--covariants of degree $L$ and $2L+1$ sets of secondary $\left(L\right)$--covariants of degree $L+1$ should be used to generate a free module $\mathcal{M}_1^{R_1}$ over the ring $R_1$ generated by the six primary invariants. However, the free module $\mathcal{M}_1^{R_1}$ is only a submodule of the module of all the $\left(L\right)$--covariant polynomials in the coordinates $\left(\mathbf{x}_1,\mathbf{x}_2,\mathbf{x}_3\right)$. The second rational function of $g_\sothree^{\left(b\right)}\left(3,L; \lambda \right)$ indicates that a second module $\mathcal{M}_2^{R_2}$ is needed, with $L\left(L-1\right)/2$ sets of $\left(L\right)$--covariants of degree $L$, $\left(L+1\right)\left(L-1\right)$ sets of $\left(L\right)$--covariants of degree $L+1$ and $L\left(L-1\right)/2$ sets of $\left(L\right)$--covariants of degree $L+2$. The module $\mathcal{M}_2^{R_2}$ should be free over a subring $R_2 \subsetneq R_1$ spanned by only five primary invariants.

We now propose a generalized integrity basis for the non--free module of $\left(2\right)$--covariants in three vectors, whose Molien generating function is:
\begin{equation}
g_\sothree^{\left(b\right)}\left(3,2; \lambda \right)
=
\frac{5\lambda^2+5\lambda^3}{(1-\lambda^2)^6}+\frac{\lambda^2+3\lambda^3+\lambda^4}{(1-\lambda^2)^5}.
\label{g2-min}
\end{equation}
The six primary invariants associated with the first rational function are chosen as the $L=0$ or $L=1$ cases, \textit{i.e.} $\left\{Q_{1,1},Q_{1,2},Q_{1,3},Q_{2,2},Q_{2,3},Q_{3,3}\right\}$. For $1\leq i\leq j\leq 3$, we define the $D_{2,M}^{i,j}$'s as:
$$
D_{2,M}^{i,j} = f_{2,M}^{(2)}\left(\mathbf{x}_i,\mathbf{x}_j\right).
$$
The $6$ $D_{2,M}^{i,j}$'s are the secondary $\left(2\right)$--covariants of degree $2$ that we are looking for. Next, we construct a similar expression by substituting $\left(x_j,y_j,z_j\right)$ with a $(L=1)$--covariant of degree $2$:
$$
\left(\begin{array}{c}
f_{2,2}^{(3)}\left(\mathbf{x}_i,\mathbf{x}_j,\mathbf{x}_k\right) \\
f_{2,1}^{(3)}\left(\mathbf{x}_i,\mathbf{x}_j,\mathbf{x}_k\right) \\
f_{2,0}^{(3)}\left(\mathbf{x}_i,\mathbf{x}_j,\mathbf{x}_k\right) \\
f_{2,-1}^{(3)}\left(\mathbf{x}_i,\mathbf{x}_j,\mathbf{x}_k\right) \\
f_{2,-2}^{(3)}\left(\mathbf{x}_i,\mathbf{x}_j,\mathbf{x}_k\right)
\end{array}\right)
=
\left( \begin{array}{c}
\frac{\sqrt{3}}{2} \left(x_i\left(y_jz_k-z_jy_k\right)-y_i\left(z_jx_k-x_jz_k\right)\right) \\
\frac{\sqrt{3}}{2}\left(x_i\left(x_jy_k-y_jx_k\right)+z_i\left(y_jz_k-z_jy_k\right)\right) \\
\frac{1}{2} \left(2z_i\left(x_jy_k-y_jx_k\right)-x_i\left(y_jz_k-z_jy_k\right)-y_i\left(z_jx_k-x_jz_k\right)\right) \\
\frac{\sqrt{3}}{2} \left(y_i\left(x_jy_k-y_jx_k\right)+z_i\left(z_jx_k-x_jz_k\right)\right) \\
\frac{\sqrt{3}}{2} \left(x_i\left(z_jx_k-x_jz_k\right)+y_i\left(y_jz_k-z_jy_k\right)\right)
\end{array} \right).
$$
However, only $8$ can be linearly independent, since $\forall \, i,j,k\in\{1,2,3\}$ and $M \in \left\{2,1,0,-1,-2\right\}$: 
$$
f_{2,M}^{(3)}\left(\mathbf{x}_i,\mathbf{x}_j,\mathbf{x}_k\right) + f_{2,M}^{(3)}\left(\mathbf{x}_i,\mathbf{x}_k,\mathbf{x}_j\right) = 0,
$$
and
$$
f_{2,M}^{(3)}\left(\mathbf{x}_i,\mathbf{x}_j,\mathbf{x}_k\right) + f_{2,M}^{(3)}\left(\mathbf{x}_j,\mathbf{x}_k,\mathbf{x}_i\right) + f_{2,M}^{(3)}\left(\mathbf{x}_k,\mathbf{x}_i,\mathbf{x}_j\right) = 0.
$$
To shorten the notation, we define $T_{2,M}^{i,j,k} = f_{2,M}^{(3)}\left(\mathbf{x}_i,\mathbf{x}_j,\mathbf{x}_k\right)$. For a given $i$ and $M$, we retain only those $T_{2,M}^{i,j,k}$ with $j<k$, and for $i,j,k$ all distinct, we decide to discard 
\begin{eqnarray}
T_{2,M}^{2,1,3}=T_{2,M}^{1,2,3}+T_{2,M}^{3,1,2}.
\label{T-linear-dep}
\end{eqnarray}
This means that we select
$$
T_{2,M}^{1,1,2}, T_{2,M}^{1,1,3}, T_{2,M}^{1,2,3}, T_{2,M}^{2,1,2}, T_{2,M}^{2,2,3}, T_{2,M}^{3,1,2}, T_{2,M}^{3,1,3}, T_{2,M}^{3,2,3},
$$
as the $8$ linearly independent $T_{2,M}^{i,j,k}$.

The three syzygies of degree $5$ are found to be, by solving linear systems of equations:
\begin{eqnarray}
Q_{2,3}T_{2,M}^{1,1,2}-Q_{2,2}T_{2,M}^{1,1,3}-Q_{1,3}T_{2,M}^{2,1,2}-Q_{1,1}T_{2,M}^{2,2,3}+Q_{1,2}(2T_{2,M}^{1,2,3}+T_{2,M}^{3,1,2}) &=& 0,
\label{siz5a} \\
Q_{1,3}T_{2,M}^{2,2,3}+Q_{3,3}T_{2,M}^{2,1,2}-Q_{1,2}T_{2,M}^{3,2,3}+Q_{2,2}T_{2,M}^{3,1,3}-Q_{2,3}(2T_{2,M}^{3,1,2}+T_{2,M}^{1,2,3}) &=& 0,
\label{siz5b} \\
Q_{1,1}T_{2,M}^{3,2,3}+Q_{2,3}T_{2,M}^{1,1,3}-Q_{1,2}T_{2,M}^{3,1,3}-Q_{3,3}T_{2,M}^{1,1,2}+Q_{1,3}(2T_{2,M}^{3,1,2}-T_{2,M}^{2,1,3}) &=& 0.
\end{eqnarray}
The last one can be rewritten by using the relation~(\ref{T-linear-dep}) as 
\begin{equation}
Q_{1,1}T_{2,M}^{3,2,3}+Q_{2,3}T_{2,M}^{1,1,3}-Q_{1,2}T_{2,M}^{3,1,3}-Q_{3,3}T_{2,M}^{1,1,2}+Q_{1,3}(T_{2,M}^{3,1,2}-T_{2,M}^{1,2,3})=0.
\label{siz5c}
\end{equation}
Similarly, the syzygy of degree $6$ is found to be:
\begin{eqnarray}
 \lefteqn{(Q_{2,3}^2-Q_{2,2}Q_{3,3})D_{2,M}^{1,1}+(Q_{1,3}^2-Q_{1,1}Q_{3,3})D_{2,M}^{2,2}+(Q_{1,2}^2-Q_{1,1}Q_{2,2})D_{2,M}^{3,3}+}
 \nonumber\\
 &2\left((Q_{1,2}Q_{3,3}-Q_{1,3}Q_{2,3})D_{2,M}^{1,2}
+(Q_{1,3}Q_{2,2}-Q_{1,2}Q_{2,3})D_{2,M}^{1,3}+(Q_{1,1}Q_{2,3}
-Q_{1,2}Q_{1,3})D_{2,M}^{2,3}\right)=0.&
\nonumber\\
\label{siz6}
\end{eqnarray}

There is some arbitrariness in the choice of the primary invariant to be removed from the ring of invariants $R_1=\mathbb{C}\left[Q_{1,1},Q_{1,2},Q_{1,3},Q_{2,2},Q_{2,3},Q_{3,3}\right]$ to define the subring $R_2 \subsetneq R_1$ of the second submodule $\mathcal{M}_2^{R_2}$. Let us remove $Q_{2,3}$, so that $R_2$ is the subring of $R_1$ generated by $\{Q_{1,1},Q_{2,2},Q_{3,3},Q_{1,2},Q_{1,3}\}$. It is then convenient to eliminate the secondary covariant $D_{2,M}^{1,1}$ of degree $2$ from $\mathcal{M}_1^{R_1}$ and to choose $D_{2,M}^{1,1}$ and covariant $Q_{2,3}\times D_{2,M}^{1,1}$ of degree $4$ as secondary covariants for $\mathcal{M}_2^{R_2}$, because relation~(\ref{siz6}) allows one to reexpress $Q_{2,3}^2 \times D_{2,M}^{1,1}$ in terms of elements of either the $\mathcal{M}_1^{R_1}$ or the $\mathcal{M}_2^{R_2}$ module, whatever choice of degree $3$ secondary covariant partitioning is made. In fact, the same is true for any product $(Q_{2,3}^n \times D_{2,M}^{1,1})$ with $n>1$ by a repeated use of the syzygy~(\ref{siz6}). It remains to select the three secondary covariants of order $3$ to exclude from $\mathcal{M}_1^{R_1}$ and to include in $\mathcal{M}_2^{R_2}$. A natural choice is $T_{2,M}^{1,1,2}$, $(2T_{2,M}^{3,1,2}+T_{2,M}^{1,2,3})$ and $T_{2,M}^{1,1,3}$, since $Q_{2,3}T_{2,M}^{1,1,2}$, $Q_{2,3}(2T_{2,M}^{3,1,2}+T_{2,M}^{1,2,3})$ and $Q_{2,3}T_{2,M}^{1,1,3}$, are easily re--expressed with terms either in $\mathcal{M}_1^{R_1}$ or in $\mathcal{M}_2^{R_2}$, by means of the syzygies (\ref{siz5a}), (\ref{siz5b}) and (\ref{siz5c}) respectively. The 5 secondary covariants of degree $2$ and the 5 secondary covariants of degree $3$ spanning  $\mathcal{M}_1^{R_1}$ can be chosen to be $D_{2,M}^{1,2}$, $D_{2,M}^{1,3}$, $D_{2,M}^{2,2}$, $D_{2,M}^{2,3}$, $D_{2,M}^{3,3}$, $T_{2,M}^{1,2,3}$, $T_{2,M}^{2,1,2}$, $T_{2,M}^{2,2,3}$, $T_{2,M}^{3,1,3}$, and $T_{2,M}^{3,2,3}$.

With this choice, the electric quadrupole moment components have the form:
\begin{eqnarray}
q_M^\mathrm{el}\left(\mathbf{x}_1,\mathbf{x}_2,\mathbf{x}_3\right) &=& P_1^\mathrm{el}\left(Q_{1,1},Q_{1,2},Q_{1,3},Q_{2,2},Q_{2,3},Q_{3,3}\right) \times D_{2,M}^{1,2}
\nonumber \\
&&+ P_2^\mathrm{el}\left(Q_{1,1},Q_{1,2},Q_{1,3},Q_{2,2},Q_{2,3},Q_{3,3}\right) \times D_{2,M}^{1,3}
\nonumber \\
&&+ P_3^\mathrm{el}\left(Q_{1,1},Q_{1,2},Q_{1,3},Q_{2,2},Q_{2,3},Q_{3,3}\right) \times D_{2,M}^{2,2}
\nonumber \\
&&+ P_4^\mathrm{el}\left(Q_{1,1},Q_{1,2},Q_{1,3},Q_{2,2},Q_{2,3},Q_{3,3}\right) \times D_{2,M}^{2,3}
\nonumber \\
&&+ P_5^\mathrm{el}\left(Q_{1,1},Q_{1,2},Q_{1,3},Q_{2,2},Q_{2,3},Q_{3,3}\right) \times D_{2,M}^{3,3}
\nonumber \\
&&+ P_6^\mathrm{el}\left(Q_{1,1},Q_{1,2},Q_{1,3},Q_{2,2},Q_{3,3}\right) \times D_{2,M}^{1,1}
\nonumber \\
&&+ P_7^\mathrm{el}\left(Q_{1,1},Q_{1,2},Q_{1,3},Q_{2,2},Q_{3,3}\right) \times Q_{2,3} D_{2,M}^{1,1}.
\label{EQMS-3vec}
\end{eqnarray}
Meanwhile, the magnetic dipole moment function is expressed as:
\begin{eqnarray}
q_M^\mathrm{mag}\left(\mathbf{x}_1,\mathbf{x}_2,\mathbf{x}_3\right) &=&
P_1^\mathrm{mag}\left(Q_{1,1},Q_{1,2},Q_{1,3},Q_{2,2},Q_{2,3},Q_{3,3}\right) \times T_{2,M}^{1,2,3}\nonumber \\
&&+ P_2^\mathrm{mag}\left(Q_{1,1},Q_{1,2},Q_{1,3},Q_{2,2},Q_{2,3},Q_{3,3}\right) \times T_{2,M}^{2,1,2}
\nonumber \\
&&+ P_3^\mathrm{mag}\left(Q_{1,1},Q_{1,2},Q_{1,3},Q_{2,2},Q_{2,3},Q_{3,3}\right) \times T_{2,M}^{2,2,3}
\nonumber \\
&&+ P_4^\mathrm{mag}\left(Q_{1,1},Q_{1,2},Q_{1,3},Q_{2,2},Q_{2,3},Q_{3,3}\right) \times T_{2,M}^{3,1,3}
\nonumber \\
&&+ P_5^\mathrm{mag}\left(Q_{1,1},Q_{1,2},Q_{1,3},Q_{2,2},Q_{2,3},Q_{3,3}\right) \times T_{2,M}^{3,2,3}
\nonumber \\
&&+ P_6^\mathrm{mag}\left(Q_{1,1},Q_{1,2},Q_{1,3},Q_{2,2},Q_{3,3}\right) \times T_{2,M}^{1,1,2}
\nonumber \\
&&+ P_7^\mathrm{mag}\left(Q_{1,1},Q_{1,2},Q_{1,3},Q_{2,2},Q_{3,3}\right) \times (2T_{2,M}^{3,1,2}+T_{2,M}^{1,2,3})
\nonumber \\
&&+ P_8^\mathrm{mag}\left(Q_{1,1},Q_{1,2},Q_{1,3},Q_{2,2},Q_{3,3}\right) \times T_{2,M}^{1,1,3}.
\label{MQMS-3vec}
\end{eqnarray}
These expansions are new and useful, for they are as compact as possible.

\subsection{\label{S3B}Other non--free modules}

The generalized integrity basis for the non--free module with three vectors and $L=2$ presented in Section~\ref{S3A} corresponds to a  decomposition of the non--free module into two free submodules. The first submodule is a module on the ring of six primary invariants while the second submodule is a module on a subring of five primary invariants only. This algebraic structure can be inferred from the Molien generating function written as a sum of two rational functions. The decomposition of a non--free module into a sum of free modules on subrings of the ring of primary invariants may have more than two terms when considering higher values of $N$. Appendix~\ref{appC} holds the Molien generating functions for four vectors amenable to a symbolic interpretation in term of an integrity basis or a generalized integrity basis. It suggests a decomposition in three rational functions whenever $L \geq 3$ (forms $g_\sothree^{(b)}\left(4,L;\lambda\right)$, $g_\sothree^{(c)}\left(4,L;\lambda\right)$ and $g_\sothree^{(d)}\left(4,L;\lambda\right)$). Appendix~\ref{appD}, for five vectors, suggests a decomposition of the Molien generating functions in four rational functions whenever $L\ge 4$ (forms $g_\sothree^{(b)}\left(5,L;\lambda\right)$ to $g_\sothree^{(g)}\left(5,L;\lambda\right)$). It is noteworthy that the suitable expression when dealing with a non--free module depends on the specific value of $L$. Table~\ref{table_1} states the Molien generating function that is to be used for a given value of $L$ as a guide in the construction of the integrity basis or generalized integrity basis.

\begin{table}
\caption{\label{table_1}Choice as a function of $L$ of the Molien generating function $g_\sothree^{\left(x\right)}\left(N,L;\lambda\right)$ for $N=3$, $4$ or $5$ vectors suitable for a symbolic interpretation in term of an integrity basis ($x=a$) or a generalized integrity basis ($x \neq a$) .}
\begin{tabular}{lll|lll|lll}
\hline\noalign{\smallskip}
\multicolumn{3}{c|}{$N=3$} & \multicolumn{3}{c|}{$N=4$} & \multicolumn{3}{c}{$N=5$} \\
$L$ & algebraic structure & $x$ & $L$ & algebraic structure & $x$ & $L$ & algebraic structure & $x$ \\
\noalign{\smallskip}\hline\noalign{\smallskip}
$0$ & ring & $a$ & $0$ & ring & $a$ & $0$ & ring & $a$ \\
$1$ & free module & $a$ & $1$ & free module & $a$ & $1-3$ & free module & $a$ \\
$2-\infty$ & non--free module & $b$ & $2$ & free module & $a$ or $b$ & $4$ & non--free module & $b$ \\
&&& $3-4$ & non--free module & $b$ & $5-6$ & non--free module & $c$ \\
&&& $5-16$ & non--free module & $c$ & $7-10$ & non--free module & $d$ \\
&&& $17-\infty$ & non--free module & $d$ & $11-14$ & non--free module & $e$ \\
&&& &&& $15-81$ & non--free module & $f$ \\
&&& &&& $82-\infty$ & non--free module & $g$ \\
\noalign{\smallskip}\hline
\end{tabular}
\end{table}

The generalized integrity bases are not unique, since, as we have seen, there is some arbitrariness in the choice of the generators of the primary invariant subrings and of the covariant generators of the modules on these subrings (although some choices are more practical than others).

\section{\label{S4}Two conjectures on the Molien generating function for the action of the orthogonal group on a set of three dimensional vectors}

\subsection{\label{S4A}First conjecture}

According to our observations described in Section~\ref{S2} and the expressions of the Molien generating function for four vectors (appendix~\ref{appC}) and five vectors (appendix~\ref{appD}), our first conjecture deals with free modules:

{\it
Conjecture 1. The generating Molien function for the action of the orthogonal group on a set of three dimensional vectors can be written as a single rational function with non--negative coefficients in the numerator for all $L$ if $N\in\left\{1,2\right\}$ and for $0 \leq L \leq N-2$ if $N \geq 3$:
$$
g_\sothree\left(N,L;\lambda \right)
=
 \frac{\sum_{n=0}^{n_{\max}\left(N,L\right)} c_{N,L}^{n}\lambda^{L+n}}{\left(1-\lambda^2\right)^{3N-3}}.
$$
This function admits a symbolic interpretation in term of an integrity basis with $3N-3$ primary polynomials of degree $2$ and $c_{N,L}^n$ secondary polynomials of degree $L+n$, $0 \leq n \leq n_{\max}\left(N,L\right)$.

In particular, the relevant covariant modules for the electric or magnetic dipole moment hypersurface corresponds to the $L=1$ irreducible representation of $\sothree$ and will always be a free module over an h.s.o.p.. 
}

The case $L=0$ corresponds to the ring of invariant polynomials $R$ under the action of the compact $\sothree$ group. 
We have seen that it is a finitely generated, free module over the algebra $\mathbb{R}[g_1, \ldots, g_m]$ for a suitable maximal set of homogeneous, algebraically independent, invariant polynomials, $\{g_1, \ldots, g_m\}$, the primary invariants \cite{HR1974,DK}. That is to say, there exists a finite set of secondary invariants $\{f_1, \ldots, f_p\}$ such that, 
$$
R = \mathbb{R}[g_1,...,g_m]f_1\oplus \cdots\oplus \mathbb{R}[g_1,...,g_m]f_p.
$$
This relationship constitutes an \textit{Hironaka decomposition} of $R$. The conjecture holds in that case as a consequence of Ref. \cite{b100}: the Molien function for invariant polynomials can be cast in the form of a single rational function. This corresponds to the fact that one can choose a h.s.o.p. with $3N-3$ primary or denominator invariant polynomials of degree $2$ and $c_{N,0}^{n}$ secondary or numerator invariant polynomials of degree $n$ with $0 \leq n \leq n_{\max}\left(N,0\right)$.

If the module of $\left(L\right)$--covariants ($L>0)$ is Cohen--Macaulay, that is, free over a suitable h.s.o.p., the situation is similar to the ring of invariant polynomials: the generating Molien function is a single rational function and the primary polynomials can be chosen to be the same as the ones for $L=0$. It is suited to an interpretation in term of integrity basis as the terms of the numerator correspond to the secondary covariants.

\subsection{\label{S4B}Second conjecture}

\subsubsection{\label{S4B1}Statement}

We formulate our second conjecture based on the results of Section~\ref{S3}. It corresponds to the case where the module of $\left(L\right)$--covariants is not free over an h.s.o.p., which is known to happen, see Ref.~\cite{VanDenBergh}.

{\it
Conjecture 2. For any number $N$ of three dimensional vectors and any final representation $L \geq N-1$ of $\sothree$, the Molien function $g_\sothree\left(N,L;\lambda\right)$ can be cast in the form of a sum of $N-1$ rational functions:

\begin{equation}
g_\sothree\left(N,L;\lambda \right)
=
\sum_{l=1}^{N-1} \frac{\mathcal{N}_{N,L}^{l}\left(\lambda\right)}{\left(1-\lambda^2\right)^{3N-2-l}},
\label{gen_molien_simple}
\end{equation}
with the $N-1$ numerator polynomials $\mathcal{N}_{N,L}^{l}\left(\lambda\right) = \sum\limits_{n=0}^{n_{\max}\left(N,L,l\right)} c_{N,L}^{l,n} \lambda^{L+n}$ having only non--negative coefficients $c_{N,L}^{l,n}$. The exponents on $\lambda$ in the numerator polynomials start at $L$, since $(L)$--covariants built from vectors are at least of total degree $L$.
} 

We already proposed to interpret this situation by introducing generalized integrity bases in our previous work on $\sotwo$ \cite{doc_03181}. The non--free module of $\left(L\right)$--covariants is then decomposed as a sum of $N-1$ submodules. The structure of the $l^\mathrm{th}$ submodule is described by the $l^\mathrm{th}$ rational function, $1\leq l\leq N-1$. Each submodule is a direct sum of free modules on rings of invariants generated by $3N-2-l$ primary invariants. In the cases we investigated, the decomposition has a stronger property: it is free by using a single set of primary invariants, that is to say, the primary invariants of the $(l+1)^\mathrm{th}$ submodule can be chosen within those of the $l^\mathrm{th}$ one. The  number of linearly independent, secondary covariants of degree $L+n$ generating those free modules is given by $c_{N,L}^{l,n}$. The set of primary invariants and secondary $\left(L\right)$--covariants associated with the $l^\mathrm{th}$ rational function is an integrity basis for the $l^\mathrm{th}$ submodule. The set of the integrity bases for all the $N-1$ submodules defines a generalized integrity basis. Our conjecture is closely related to the conjecture 5.1 of Stanley in Ref. \cite{Stanley82}, which has a more general scope.

\subsubsection{\label{S4B2}Heuristic for the rewriting of the Molien generating function in a form suitable for its symbolic interpretation in term of a generalized integrity basis}

It is desirable to have for a given pair $(N,L)$ an heuristic that algorithmically determines the numerators in the sum of expression~(\ref{gen_molien_simple}) by starting with the Molien generating function written as a single rational function with numerator $\mathcal{N}_\sothree^{\left(a\right)}\left(N,L;\lambda\right)$. If for the value of $L$ considered, the numerator $\mathcal{N}_\sothree^{\left(a\right)}\left(N,L;\lambda\right)$ has only non--negative coefficients, then the form is suitable for an interpretation in term of integrity basis and we can start the construction of the integrity basis. Otherwise we perform a modified polynomial long division of the numerator by $1-\lambda^2$. The modification is in the halt criterion: the division process is stopped not when the degree of the remainder $r_1$ is less than $2$ as one would do in the usual polynomial division, but when the coefficients of the remainder $r_1$ become all non--negative for the value of $L$ considered. Then the numerator is rewritten as:
$$
\mathcal{N}_\sothree^{\left(a\right)}\left(N,L;\lambda\right) = r_1 + \left(1-\lambda^2\right) q_1.
$$
The remainder $r_1$ serves as the numerator of the first ($l=1$) rational function in the right--hand side of~(\ref{gen_molien_simple}). If all the coefficients in the quotient $q_1$ are non--negative for the value of $L$ considered, the division stops and the final form of the Molien function is then:
$$
g_\sothree^{\left(b\right)}\left(N,L;\lambda\right) = \frac{r_1}{\left(1-\lambda^2\right)^{3N-3}}
  +\frac{q_1}{\left(1-\lambda^2\right)^{3N-4}}.
$$
  If at least one coefficient in the quotient $q_1$ is negative, the modified polynomial division procedure is applied to $q_1$. The quotient $q_1$ is decomposed as $q_1 = r_2  + \left(1-\lambda^2\right) q_2$ and the new remainder $r_2$ with non--negative coefficients will constitute the numerator of the second ($l=2$) rational function. The modified polynomial division is repeated until the quotient has also only non--negative coefficients.

For example, the single rational function (\ref{g_three_vectors_1}) with its numerator~(\ref{g_three_vectors_1_}) is suitable for $N=3$ and $L \in \left\{0,1\right\}$. For $L\geq 2$, we obtain by the modified polynomial division of the numerator~(\ref{g_three_vectors_1_}):
\begin{eqnarray}
\lefteqn{
\frac{(L+2)(L+1)}{2}\lambda^L +(L+2)L\lambda^{L+1} 
-(L+1)(L-1)\lambda^{L+3} -\frac{L(L-1)}{2}\lambda^{L+4}
} \nonumber \\
&&=(1-\lambda^2)\left(\frac{L(L-1)}{2}\lambda^{L}+(L+1)(L-1)\lambda^{L+1}+\frac{L(L-1)}{2}\lambda^{L+2} \right)+(2L+1)\lambda^L +(2L+1)\lambda^{L+1}.
\nonumber \\
\label{algo-3-vec}
\end{eqnarray}
The remainder $(2L+1)\lambda^L +(2L+1)\lambda^{L+1}$ has only non--negative coefficients and is the numerator of the first rational function we are seeking for. All the coefficients of the quotient $\frac{L(L-1)}{2}\lambda^{L}+(L+1)(L-1)\lambda^{L+1}+\frac{L(L-1)}{2}\lambda^{L+2}$ are non--negative for $L\geq 2$. This quotient is the numerator of the second fraction. This procedure gives again the expression (\ref{g_three_vectors_2}) of the Molien function written as a sum of two rational functions with non--negative coefficients in the numerators.

For $N=4$, all the coefficients of the numerator $\mathcal{N}_\sothree^{\left(a\right)}\left(4,L;\lambda\right)$ of (\ref{num_four_vectors_1}) are non--negative for $L<3$. For $L\geq 3$, we obtain by the modified polynomial division of the initial numerator:
\begin{scriptsize}
\begin{eqnarray}
\lefteqn{
\frac{(L+3)(L+2)(L+1)}{6} \lambda^L 
+ \frac{(L+3)(L+2)L}{2} \lambda^{L+1}
+ \frac{(L+3)(L+2)(L+1)}{6} \lambda^{L+2}
- \frac{(L+3)(L-2)(5L+4)}{6} \lambda^{L+3}
}
\nonumber \\
\lefteqn{
- \frac{(L+3)(L-2)(5L+1)}{6} \lambda^{L+4}
+ \frac{L(L-1)(L-2)}{6} \lambda^{L+5}
+ \frac{(L+1)(L-1)(L-2)}{2} \lambda^{L+6}
+ \frac{L(L-1)(L-2)}{6} \lambda^{L+7}
}
\nonumber \\
&&=(1-\lambda^2)\left( \frac{(L-2)(L^2+7L+4)}{2} \lambda^{L+1}+\frac{(L-2)(L^2+8L+3)}{3} \lambda^{L+2}-\frac{L(L-1)(L-2)}{3}\lambda^{L+3}- \frac{(L+1)(L-1)(L-2)}{2} \lambda^{L+4}  \right.
\nonumber \\
&&
 \left. -\frac{L(L-1)(L-2)}{6}\lambda^{L+5}  \right)
+\frac{(L+3)(L+2)(L+1)}{6} \lambda^L 
+ 4(2L+1) \lambda^{L+1}
- \frac{(L^3+6L^2-37L-18)}{6} \lambda^{L+2}
\label{algo-4-vec1}
\\
&&=(1-\lambda^2)\left( \frac{(L^3+6L^2-37L-18)}{6} \lambda^{L}+\frac{(L-2)(L^2+7L+4)}{2} \lambda^{L+1}+\frac{(L-2)(L^2+8L+3)}{3} \lambda^{L+2}-\frac{L(L-1)(L-2)}{3}\lambda^{L+3}  \right.
\nonumber \\
&&
 \left. - \frac{(L+1)(L-1)(L-2)}{2} \lambda^{L+4}-\frac{L(L-1)(L-2)}{6}\lambda^{L+5}  \right)
+ 4(2L+1) \lambda^L 
+ 4(2L+1) \lambda^{L+1}.
\label{algo-4-vec2}
\end{eqnarray}
\end{scriptsize}
For $L \in \left\{3,4\right\}$, the division stops at~(\ref{algo-4-vec1}), where all coefficients of the remainder are non--negative, whereas for larger values of $L$ the coefficient $-\left(L^3+6L^2-37L-18\right)/6$ of $\lambda^{L+2}$ is negative in the remainder of~(\ref{algo-4-vec1}) and one must stop a step later at~(\ref{algo-4-vec2}). Let us consider first the cases $L=3$ and $L=4$. The quotient has negative coefficients and must be divided again by $(1-\lambda^2)$,\\
\begin{scriptsize}
\begin{eqnarray}
\lefteqn{
\frac{(L-2)(L^2+7L+4)}{2} \lambda^{L+1}+\frac{(L-2)(L^2+8L+3)}{3} \lambda^{L+2}-\frac{L(L-1)(L-2)}{3}\lambda^{L+3}- \frac{(L+1)(L-1)(L-2)}{2} \lambda^{L+4}   -\frac{L(L-1)(L-2)}{6}\lambda^{L+5}
} \nonumber \\
&&=(1-\lambda^2)\left(\frac{L(L-1)(L-2)}{2}\lambda^{L+1}+\frac{(L+1)(L-1)(L-2)}{2} \lambda^{L+2}+\frac{L(L-1)(L-2)}{6}\lambda^{L+3}  \right)+2(L-2)(2L+1) \lambda^{L+1}
\nonumber \\
&&
-\frac{(L-2)(L^2-16L-9)}{6} \lambda^{L+2}.
\label{algo-4-vec3}
\end{eqnarray}
\end{scriptsize}
Both remainder and quotient have non--negative coefficients for $L\in\left\{3,4\right\}$, so we stop here and obtain again the numerators of the second and third fractions of $g_\sothree^{\left(b\right)}\left(4,L;\lambda\right)$ in~(\ref{gen_four_vectors_2}).

For $L>4$, we have to divide the quotient of~(\ref{algo-4-vec2}) by $(1-\lambda^2)$,
\begin{scriptsize}
\begin{eqnarray}
\lefteqn{
\frac{(L^3+6L^2-37L-18)}{6} \lambda^{L}+\frac{(L-2)(L^2+7L+4)}{2} \lambda^{L+1}+\frac{(L-2)(L^2+8L+3)}{3} \lambda^{L+2}-\frac{L(L-1)(L-2)}{3}\lambda^{L+3}
} \nonumber \\
&&
- \frac{(L+1)(L-1)(L-2)}{2} \lambda^{L+4}-\frac{L(L-1)(L-2)}{6}\lambda^{L+5}
\nonumber \\
&& =(1-\lambda^2)\left(\frac{L(L-1)(L-2)}{2}\lambda^{L+1}+\frac{(L+1)(L-1)(L-2)}{2} \lambda^{L+2}+\frac{L(L-1)(L-2)}{6}\lambda^{L+3}  \right)+ \frac{(L^3+6L^2-37L-18)}{6} \lambda^{L}
\nonumber \\
&& +2(L-2)(2L+1) \lambda^{L+1}-\frac{(L-2)(L^2-16L-9)}{6} \lambda^{L+2}
\label{algo-4-vec4}
\\
&& =(1-\lambda^2)\left(\frac{(L-2)(L^2-16L-9)}{6} \lambda^{L}+\frac{L(L-1)(L-2)}{2}\lambda^{L+1}+\frac{(L+1)(L-1)(L-2)}{2} \lambda^{L+2}+\frac{L(L-1)(L-2)}{6}\lambda^{L+3}  \right)
\nonumber
\\
&& + 2(L-3)(2L+1) \lambda^{L}+2(L-2)(2L+1) \lambda^{L+1}.
\label{algo-4-vec5}
\end{eqnarray}
\end{scriptsize}
For $5\leq L\leq 16$  the division stops at~(\ref{algo-4-vec4}), when all coefficients of the remainder are non--negative. For larger values of $L$, the coefficient of $\lambda^{L+2}$ is negative in the remainder of~(\ref{algo-4-vec4}) and the division must continue up to~(\ref{algo-4-vec5}). In both cases, the quotient has only non--negative coefficients, so the remainders and quotients of Eqs.~(\ref{algo-4-vec4}) and~(\ref{algo-4-vec5}) gives the numerators of the second and third fractions of $g_\sothree^{\left(c\right)}\left(4,L;\lambda\right)$ and $g_\sothree^{\left(d\right)}\left(4,L;\lambda\right)$ in Eqs.~(\ref{gen_four_vectors_3})  and~(\ref{gen_four_vectors_4}).

\subsubsection{\label{S4B3}Construction of a generalized integrity basis}

The heuristic presented in Section~\ref{S4B2} provides an essential help in the construction of the generalized integrity basis of a non--free module. To the sum in the right--hand side of~(\ref{gen_molien_simple}) corresponds a decomposition of the non--free module in a sum of $N-1$ submodules $\mathcal{M}_l^{R_l}$ over a ring $R_l$:
\begin{equation}
\label{decomposition}
\mathcal{M} = \oplus_{l=1}^{N-1} \mathcal{M}_l^{R_l}.
\end{equation}
The number and degrees of the secondary covariants which are the basis of each free submodule $\mathcal{M}_l^{R_l}$ in the decomposition~(\ref{decomposition}) are given by the successive remainders of the divisions by $\left(1-\lambda^2\right)$ and the last quotient with only positive coefficients. To find the successive sets of syzygies, it is enough to solve linear systems in the primary invariant ring and its successively selected subrings. The numbers of independent syzygies to write and their degrees are given by the negative coefficients in the expressions of the numerator $\mathcal{N}_\sothree^{\left(a\right)}\left(N,L,\lambda\right)$ and its successive quotients by $(1-\lambda^2)$ appearing while following the heuristic of Section~\ref{S4B2}. For example, for $N=4$, $L=3$, Eq.~(\ref{num_four_vectors_1}) tells us that there will be $19$ syzygies of degree $6$ and $16$ of degree $7$ to be used in order to select  $R_2$ and the $20$ secondary covariants of degree $3$, $28$ secondary covariants of degree $4$ and $8$ secondary covariants of degree $5$ of $\mathcal{M}_1^{R_1}$ according to the remainder in~(\ref{algo-4-vec1}). This first set of syzygies is to be obtained by solving a linear system in $R_1$. Then, the quotient of~(\ref{algo-4-vec1}) tells us that a second set of $2$ syzygies of degree  $6$, $4$ of degree $7$ and $1$ of degree $8$ is to be obtained by solving linear systems in the subring $R_2$, previously selected. Finally, the remainder in~(\ref{algo-4-vec3}), tells us that there are $14$ secondary covariants of degree $4$ and $8$ of degree $5$ to be chosen for $\mathcal{M}_2^{R_2}$, and the quotient in~(\ref{algo-4-vec3}) that, once  $R_3$ has been selected, there will be $3$ secondary covariants of degree $4$, $4$ of degree $5$ and $1$ of degree $6$ to be found for $\mathcal{M}_3^{R_3}$.
 
Since the $\mathcal{S}_N$ permutation group action on vector indices preserves partial degrees, one can take advantage of the Molien functions with distinguished representation arguments to obtain information about the partial degrees of the variables in these syzygies. This reduces significantly the size of the linear systems to be solved. For example, in the $N=3$ and $L=2$ case, the syzygy of order $6$ is found to be of partial degrees $n_1=n_2=n_3=2$ from the last term in~(\ref{N_three_vectors_1_L2}).  We however have not systematically reported such detailed expressions for $N>2$ because the closed formulas are not polynomial. For example, for $N=3$, we have obtained the following numerator
\begin{eqnarray}
\lefteqn{\mathcal{N}_\sothree^{\left(a\right)}\left(3,L;\lambda_1,\lambda_2,\lambda_3 \right)}
\nonumber \\
&=&\frac{1}{\left(\lambda_1-\lambda_2\right)\left(\lambda_1-\lambda_3\right) 
\left(\lambda_2-\lambda_3\right)}\left[\lambda_2\lambda_3 \left(-\lambda_3^{1+L}-\lambda_2 \lambda_3^{1+L}+\lambda_2^{1+L} \left(1+\lambda_3\right)\right)\right.
\nonumber\\
&&+\lambda_1^3\lambda_2 \lambda_3 \left(-\lambda_3^{1+L}-\lambda_2 \lambda_3^{1+L}+\lambda_2^{1+L} \left(1+\lambda_3\right)\right)
\nonumber \\
&&+\lambda_1^{2+L} \left(-\lambda_3 \left(1+\lambda_3\right)-
\lambda_2^3 \lambda_3 
\left(1+\lambda_3\right)+\lambda_2 \left(1+\lambda_3^3\right)+\lambda_2^2 \left(1+\lambda_3^3\right)\right)
\nonumber \\
&&\left.
+\lambda_1 \left(\lambda_3^{2+L}+\lambda_2^3 \lambda_3^{2+L}-\lambda_2^{2+L} \left(1+\lambda_3^3\right)\right)+\lambda_1^2
\left(\lambda_3^{2+L}+\lambda_2^3 \lambda_3^{2+L}-\lambda_2^{2+L}\left(1+\lambda_3^3\right)\right)\right]
\nonumber \\
\label{N_three_vectors_1_L} 
\end{eqnarray}
which we only managed to simplify into interpretable polynomial expressions such as~(\ref{N_three_vectors_1_L2}), when $L$ takes specific values. Such interpretable detailed expressions are useful to provide the partial degrees of the secondary covariants of the modules, after division by the factor in the denominator corresponding to the primary invariant excluded from the next subring in the decomposition: for example, if the invariant $Q_{2,3}$ is excluded one has to divide by $(1-\lambda_2 \lambda_3)$.

\section{\label{S5}Conclusion}

Expressions of the Molien generating function for the action of the $\sothree$ group on a set of three dimensional vectors have been given. Such expressions are useful guides for the construction of invariant and covariant (possibly generalized) integrity bases. The extension to the $\othree$ group is direct. When the module of covariants is non--free over an h.s.o.p., an heuristic has been proposed to transform the Molien function written as a single rational function into a form amenable to a symbolic interpretation in term of a generalized integrity basis. The same heuristic can guide step by step the construction of such a generalized integrity basis. Within this approach, the non--free module is decomposed as a direct sum of submodules. Each submodule is associated with a rational function indicating the numbers of primary invariants and secondary covariants, the former decreasing in the consecutive rational functions.

The case of the invariant ("$L=0$ covariants") module is always free \cite{stanleycov}. In quantum physics, an integrity basis can be useful to express $\sothree$--totally invariant observables such as the so--called potential energy (hyper)surface (PES) in quantum chemistry \cite{a1776}, when it is not possible or appropriate to separate out rotational from internal coordinates. This is the case of very floppy molecules such as CH$_5^+$, in order to use its full symmetry group \cite{Schmiedt15}. Note that PES are actually $\othree$--totally invariant, but for polyatomic molecules of four atoms and more, it is more practical to parametrize them by using coordinates that are $\sothree$-- but not $\othree$--invariants, such as dihedral angles.

We have conjectured that the module of $L=1$ covariants will always be free as well. We have provided explicit integrity bases up to $N=3$. In quantum physics, these integrity bases can be useful to express observables such as a dipole moment hypersurface, used in theoretical spectroscopy to calculate dipolar transition intensities.  
If there is a finite group action on the vector variables in addition to the $\sothree$ action, we can take advantage of it too as was argued in Ref.~\cite{a1776} for the particular case of invariants. For example, for $N=2$, related to the case of a triatomic molecule $\mathrm{ABC}$, if the origin of the two vectors is $\mathrm{A}$ and if atoms $\mathrm{B}$ and $\mathrm{C}$ are equivalent, then the action of the permutation group $\mathcal{S}_2$ on the vectors $\overrightarrow{\mathrm{AB}} = x_1 \, \hat{e}_x + y_1 \, \hat{e}_y + z_1 \, \hat{e}_z$ and $\overrightarrow{\mathrm{AC}} = x_2 \, \hat{e}_x + y_2 \, \hat{e}_y + z_2\, \hat{e}_z$ can be exploited to simplify the expression of the physical observables such as the dipole moment surface functions. In this particular case, we deduce for example that the polynomials of~(\ref{EDMS}) must satisfy $P_1^\mathrm{el}\left(Q_{1,1},Q_{2,2},Q_{1,2}\right)=P_2^\mathrm{el}\left(Q_{2,2},Q_{1,1},Q_{1,2}\right)$.

The $(L=2)$--case is useful in molecular physics and chemistry to expand quadrupole moments. Since the first detection of the electric quadrupole spectrum in molecular hydrogen \cite{herzberg}, most of the experimental data dealt with homonuclear diatomics until recently. The combination of accurate variational calculations and sensitive spectroscopic experiments make now possible the detection of electric quadrupole infrared transitions in triatomic molecules such as water \cite{camparguePRR,camparguePCCP} and carbon dioxide \cite{yachmenev}. One of the main achievement of this work has been to provide compact expansions for quadrupole moment hypersurfaces, which take full advantage of the structure of the non--free covariant module through the use of a generalized integrity basis.
However, the elimination of some primary invariants breaks the permutational symmetry of the subrings used for the modules $\mathcal{M}_2^{R^2}$, $\mathcal{M}_3^{R^3}$, $\ldots$, of the decomposition of the module of covariants. So, unfortunately in such a case, further permutational symmetry adaptation may be unpractical in general.


\begin{acknowledgements}
The authors thank Boris I. Zhilinski\'{\i} for fruitful discussions. PCC and FP are indebted to Prof. C. Procesi for noticing a mistake in the description of the Hironaka decomposition in Ref.[13].
The statement has been amended in the present article. 
\end{acknowledgements}

\section*{Funding}
Financial support for the project \textit{Application de la Th\'eorie des Invariants \`a la Physique Mol\'eculaire} via a CNRS grant \textit{Projet Exploratoire Premier Soutien} (PEPS) \textit{Physique Th\'eorique et Interfaces} (PTI) is acknowledged.

\section*{Conflicts of interest/Competing interests}
The authors have no conflicts of interest to declare that are relevant to the content of this article.

\section*{Availability of data and material}
Not applicable.

\section*{Code availability}
Not applicable.

\section*{Authors' contributions}
All authors contributed to the study conception. The investigation was performed by Guillaume Dhont and Patrick Cassam--Chena{\"\i}. The first draft of the manuscript was written by Guillaume Dhont and Patrick Cassam--Chena{\"\i} and all authors commented on previous versions of the manuscript. All authors read and approved the final manuscript.

\appendix

\section{\label{appA}Determination of the Molien generating function}

\subsection{Expression via an integral}

In an active transformation, a rotation of angle $\omega$ around a rotation axis whose direction is given by the unit vector $\hat{n}=n_x \, \hat{e}_x + n_y \, \hat{e}_y + n_z \, \hat{e}_z$ transforms the vector $\vec{r} = x \, \hat{e}_x + y \, \hat{e}_y + z \, \hat{e}_z$ into the vector $\vec{r^\prime} = x^\prime \, \hat{e}_x + y^\prime \, \hat{e}_y + z^\prime \, \hat{e}_z$. The initial and final coordinates are related by \cite{BL}:
$$
\left( \begin{array}{c} x^\prime \\ y^\prime \\ z^\prime \end{array} \right)
=
R\left(\omega,\hat{n}\right)
\left( \begin{array}{c} x \\ y \\ z \end{array} \right),
\quad
R\left(\omega,\hat{n}\right) = 1_3 + N \sin \omega + N^2 \left(1-\cos \omega\right),
$$
where $1_3$ is the $3\times 3$ identity matrix and $N$ is the $3\times 3$ matrix defined as:
$$
N = \left( \begin{array}{ccc}
0 & -n_z & n_y \\
n_z & 0 & -n_x \\
-n_y & n_x & 0
\end{array}
\right).
$$
Using spherical coordinates to give the orientation of the rotation axis $\hat{n}$ gives $n_x = \sin \theta \cos \varphi$, $n_y = \sin \theta \sin \varphi$ and $n_z = \cos \theta$. The $3N\times3N$ block matrix representation $D\left(\omega,\hat{n}\right)$ of the rotation operation is:
$$
D\left(\omega,\hat{n}\right)
=\left( \begin{array}{cccc}
R\left(\omega,\hat{n}\right) & 0 & \cdots & 0 \\
0 & R\left(\omega,\hat{n}\right) & \cdots & 0 \\
\vdots & \vdots & \ddots & 0 \\
0 & 0 & \cdots & R\left(\omega,\hat{n}\right)
\end{array} \right),
$$
and one easily finds that
$$
\det \left(1_{3N}-\lambda D\left(\omega,\hat{n}\right)\right) =
\left[(1-\lambda)(1-2\lambda\cos \omega+\lambda^2)\right]^N.
$$
The Molien function~(\ref{molien1}) then reduces to
\begin{equation}
\label{genfun_01}
g_\sothree\left( N,L;\lambda \right)
= \frac{2}{\pi}
\frac{1}{(1-\lambda)^N}
\int_{0}^{\pi}
\frac{\sin \left[\left(2L+1\right)\frac{\omega}{2}\right]
\sin \frac{\omega}{2}}{(1-2\lambda\cos \omega+\lambda^2)^N}
d\omega.
\end{equation}
The integral in~(\ref{genfun_01}) can be evaluated by the use of the product--to--sum identity $\sin \left[\left(2L+1\right)\omega/2\right] \sin \left(\omega/2\right) = \left\{ \cos \left(L \omega\right)-\cos \left[\left(L+1\right) \omega\right] \right\}/2$ and the tabulated formula of Ref.~\cite{equation_3_616_7_gradshteyn_ryzhik}:
\begin{eqnarray}
\lefteqn{\int_0^\pi \frac{\cos nx \, dx}{\left(1-2a\cos x+a^2\right)^m}}
\nonumber\\
&=& 
\frac{a^{2m+n-2}\pi}{\left(1-a^2\right)^{2m-1}}
\sum_{k=0}^{m-1}
\left( \begin{array}{c} m+n-1 \\ k \end{array} \right)
\left( \begin{array}{c} 2m-k-2 \\ m-1 \end{array} \right)
\left(\frac{1-a^2}{a^2}\right)^k,
\quad a^2<1.
\label{gradshteyn}
\end{eqnarray}

\subsection{\label{fun_gen_by_coupling}Recursive formula for two or more vectors}

The Molien generating function for $N \geq 2$ can be determined from formula~(\ref{genfun_01}). There is however an other approach, which gives more insight in the forthcoming construction of the integrity bases. According to the triangular conditions of the theory of angular momentum, the coupling of a set of $\left(L_1\right)$--covariants with a set of $\left(L_2\right)$--covariants generates one set of $\left(L\right)$--covariants, with $\left|L_1-L_2\right| \leq L \leq L_1+L_2$. Correspondingly the Molien generating function for $N_1+N_2$ vectors and final representation $(L)$ can be computed from the Molien generating functions for $N_1$ and $N_2$ vectors, see equation~(43) of Ref.~\cite{doc_00628}:

\begin{equation}
\label{formule_article_a0628}
g_\sothree\left(N_1+N_2,L;\lambda_1,\lambda_2\right) =
\sum_{L_1=0}^\infty \sum_{L_2=0}^\infty \Delta\left(L_1,L_2,L\right) g_\sothree\left(N_1,L_1;\lambda_1\right) g_\sothree\left(N_2,L_2;\lambda_2\right),
\end{equation}
where $\Delta\left(L_1,L_2,L\right)=1$ if $\left|L_1-L_2\right| \leq L \leq L_1+L_2$ and $0$ otherwise.

\section{\label{realSH}Complex and real solid harmonics}

The complex solid harmonics are homogeneous polynomials of degree $l$ in $x$, $y$, and $z$. Their expression is given by \cite{BL_3.153}:
$$
\mathcal{Y}_{l,m}\left(\mathbf{x}\right)
=
\sqrt{\frac{2l+1}{4\pi} \left(l+m\right)! \left(l-m\right)!}
\sum_k \frac{\left(-x-iy\right)^{k+m} \left(x-iy\right)^k z^{l-2k-m}}{2^{2k+m} \left(k+m\right)! k! \left(l-m-2k\right)!},
$$
where $\mathbf{x}$ is the triplet of coordinates $\left(x,y,z\right)$.
The complex solid harmonics are complex--valued functions for $m \neq 0$ and satisfy the property
$$
\mathcal{Y}_{l,m}\left(\mathbf{x}\right)^* = \left(-1\right)^m \mathcal{Y}_{l,-m}\left(\mathbf{x}\right).
$$
Biedenharn and Louck define in Ref.~\cite{BL} the real solid harmonics $\bar{\mathcal{Y}}_{l,m}$ through the linear combinations of complex solid harmonics \cite{BL_6.183}:
\begin{eqnarray}
\bar{\mathcal{Y}}_{l,m}\left(\mathbf{x}\right) &=&
-\frac{1}{\sqrt{2}} \left[ \mathcal{Y}_{l,m}\left(\mathbf{x}\right) + \left(-1\right)^m \mathcal{Y}_{l,-m}\left(\mathbf{x}\right) \right],
\nonumber \\
\bar{\mathcal{Y}}_{l,0}\left(\mathbf{x}\right) &=& \mathcal{Y}_{l,0}\left(\mathbf{x}\right),
\nonumber \\
\bar{\mathcal{Y}}_{l,-m}\left(\mathbf{x}\right) &=&
\frac{i}{\sqrt{2}} \left[ \mathcal{Y}_{l,m}\left(\mathbf{x}\right) - \left(-1\right)^m \mathcal{Y}_{l,-m}\left(\mathbf{x}\right) \right],
\nonumber
\end{eqnarray}
with $m \in \left\{l,l-1,\cdots,1\right\}$.

Other linear combinations of complex solid harmonics exist in the litterature. Steinborn \cite{steinborn} or Blanco \textit{et al.} \cite{BFB} use:
\begin{eqnarray}
\tilde{\mathcal{Y}}_{l,m}\left(\mathbf{x}\right) &=&
\frac{\left(-1\right)^m}{\sqrt{2}} \left[ \mathcal{Y}_{l,m}\left(\mathbf{x}\right) + \left(-1\right)^m \mathcal{Y}_{l,-m}\left(\mathbf{x}\right) \right],
\nonumber \\
\tilde{\mathcal{Y}}_{l,0}\left(\mathbf{x}\right) &=& \mathcal{Y}_{l,0}\left(\mathbf{x}\right),
\nonumber \\
\tilde{\mathcal{Y}}_{l,-m}\left(\mathbf{x}\right) &=&
\frac{\left(-1\right)^m}{i\sqrt{2}} \left[ \mathcal{Y}_{l,m}\left(\mathbf{x}\right) - \left(-1\right)^m \mathcal{Y}_{l,-m}\left(\mathbf{x}\right) \right],
\nonumber
\end{eqnarray}
with $m \in \left\{l,l-1,\cdots,1\right\}$. The two definitions of real solid harmonics are identical except for $m$  even where the two definitions give expressions with opposite sign. In the main text, we use $\tilde{\mathcal{Y}}_{l,m}$ as real solid harmonics. Their expression for up to $L=3$ are given below:

$$
\sqrt{4\pi} \tilde{\mathcal{Y}}_{0,0}\left(\mathbf{x}\right) = 1
$$

$$
\sqrt{\frac{4\pi}{3}} 
\left( \begin{array}{c}
\tilde{\mathcal{Y}}_{1,1}\left(\mathbf{x}\right) \\
\tilde{\mathcal{Y}}_{1,0}\left(\mathbf{x}\right) \\
\tilde{\mathcal{Y}}_{1,-1}\left(\mathbf{x}\right)
\end{array} \right)
= \left( \begin{array}{c} x \\ z \\ y \end{array} \right)
$$

$$
\sqrt{\frac{4\pi}{5}}
\left( \begin{array}{c}
\tilde{\mathcal{Y}}_{2,2}\left(\mathbf{x}\right) \\
\tilde{\mathcal{Y}}_{2,1}\left(\mathbf{x}\right) \\
\tilde{\mathcal{Y}}_{2,0}\left(\mathbf{x}\right) \\
\tilde{\mathcal{Y}}_{2,-1}\left(\mathbf{x}\right) \\
\tilde{\mathcal{Y}}_{2,-2}\left(\mathbf{x}\right)
\end{array} \right)
= \left( \begin{array}{c}
\frac{\sqrt{3}}{2}\left(x^2 - y^2\right) \\
\sqrt{3}xz \\
\frac{1}{2}\left(2z^2-x^2-y^2\right) \\
\sqrt{3}yz \\
\sqrt{3}xy
\end{array} \right)
$$

$$
\sqrt{\frac{4\pi}{7}}
\left( \begin{array}{c}
\tilde{\mathcal{Y}}_{3,3}\left(\mathbf{x}\right) \\
\tilde{\mathcal{Y}}_{3,2}\left(\mathbf{x}\right) \\
\tilde{\mathcal{Y}}_{3,1}\left(\mathbf{x}\right) \\
\tilde{\mathcal{Y}}_{3,0}\left(\mathbf{x}\right) \\
\tilde{\mathcal{Y}}_{3,-1}\left(\mathbf{x}\right) \\
\tilde{\mathcal{Y}}_{3,-2}\left(\mathbf{x}\right) \\
\tilde{\mathcal{Y}}_{3,-3}\left(\mathbf{x}\right)
\end{array} \right)
= \left( \begin{array}{c}
\frac{\sqrt{10}}{4} x\left(x^2-3y^2\right) \\
\frac{\sqrt{15}}{2}z\left(x^2-y^2\right) \\
\frac{\sqrt{6}}{4} \left(4z^2-x^2-y^2\right)x \\
\frac{1}{2}z\left(2z^2-3x^2-3y^2\right) \\
\frac{\sqrt{6}}{4}\left(4z^2-x^2-y^2\right)y \\
\sqrt{15}xyz \\
\frac{\sqrt{10}}{4} y\left(3x^2-y^2\right)
\end{array} \right)
$$

\section{\label{appC}Molien generating function for four vectors}

The Molien generating function for four vectors can be written as a single rational function:
\begin{equation}
g_\sothree^{\left(a\right)}\left(4,L;\lambda\right) =
\frac{\mathcal{N}_\sothree^{\left(a\right)}\left(4,L;\lambda\right)}{( 1-\lambda^2 )^9},
\label{gen_four_vectors_1}
\end{equation}
with its numerator equal to:
\begin{eqnarray}
\mathcal{N}_\sothree^{\left(a\right)}\left(4,L;\lambda \right) &=&
\frac{(L+3)(L+2)(L+1)}{6} \lambda^L 
+ \frac{(L+3)(L+2)L}{2} \lambda^{L+1}
\nonumber \\
&&
+ \frac{(L+3)(L+2)(L+1)}{6} \lambda^{L+2}
- \frac{(L+3)(L-2)(5L+4)}{6} \lambda^{L+3}
\nonumber \\
&&
- \frac{(L+3)(L-2)(5L+1)}{6} \lambda^{L+4}
+ \frac{L(L-1)(L-2)}{6} \lambda^{L+5}
\nonumber \\
&&
+ \frac{(L+1)(L-1)(L-2)}{2} \lambda^{L+6}
+ \frac{L(L-1)(L-2)}{6} \lambda^{L+7}.
\label{num_four_vectors_1}
\end{eqnarray}

The Molien generating function~(\ref{gen_four_vectors_1}) has non--negative coefficients in its numerator for $L\in \left\{0,1,2\right\}$. Negative coefficients appear for $L\geq 3$. However, the Molien function can be rewritten as Eq.~(\ref{gen_four_vectors_2}), which has only non--negative coefficients in the numerator for $L\in\left\{2,3,4\right\}$,
\begin{eqnarray}
\lefteqn{g_\sothree^{\left(b\right)}\left(4,L;\lambda\right)}
\nonumber \\
&=&
\frac{\frac{(L+3)(L+2)(L+1)}{6} \lambda^L 
  + 4(2L+1) \lambda^{L+1}
  + (-\frac{1}{6}L^3-L^2+\frac{37}{6}L+3) \lambda^{L+2}}{(1-\lambda^2)^9}
\nonumber \\
&&
+\frac{2(L-2)(2L+1) \lambda^{L+1}
  - \frac{(L-2)(L^2-16L-9)}{6} \lambda^{L+2}}{(1-\lambda^2)^8}
\nonumber \\
&&
+\frac{\frac{L(L-1)(L-2)}{2} \lambda^{L+1}
  + \frac{(L+1)(L-1)(L-2)}{2} \lambda^{L+2}
  + \frac{L(L-1)(L-2)}{6} \lambda^{L+3}}{(1-\lambda^2)^7},
\label{gen_four_vectors_2}
\end{eqnarray}
as~(\ref{gen_four_vectors_3}), which has only non--negative coefficients in the numerator for $L$ between $5$ and $16$,
\begin{eqnarray}
\lefteqn{g_\sothree^{\left(c\right)}\left(4,L;\lambda\right)}
\nonumber \\
&=&
\frac{4(2L+1) \lambda^L 
  + 4(2L+1) \lambda^{L+1}}{(1-\lambda^2)^9}
\nonumber \\
&&
+\frac{(\frac{1}{6}L^3+L^2-\frac{37}{6}L-3) \lambda^L 
  + 2(L-2)(2L+1) \lambda^{L+1}
  - \frac{(L-2)(L^2-16L-9)}{6} \lambda^{L+2}}{(1-\lambda^2)^8}
\nonumber \\
&&
+\frac{\frac{L(L-1)(L-2)}{2} \lambda^{L+1}
  + \frac{(L+1)(L-1)(L-2)}{2} \lambda^{L+2}
  + \frac{L(L-1)(L-2)}{6} \lambda^{L+3}}{(1-\lambda^2)^7},
\label{gen_four_vectors_3}
\end{eqnarray}
and as Eq.~(\ref{gen_four_vectors_4}), which has only non--negative coefficients in the numerator for $L\geq 17$,
\begin{eqnarray}
\lefteqn{g_\sothree^{\left(d\right)}\left(4,L;\lambda\right)}
\nonumber \\
&=&
\frac{4(2L+1) \lambda^L 
  + 4(2L+1) \lambda^{L+1}}{(1-\lambda^2)^9}
\nonumber \\
&&
+\frac{2(L-3)(2L+1) \lambda^L
  + 2(L-2)(2L+1) \lambda^{L+1}}{(1-\lambda^2)^8}
\nonumber \\
&&
+\frac{\frac{(L-2)(L^2-16L-9)}{6} \lambda^L
  + \frac{L(L-1)(L-2)}{2} \lambda^{L+1}
  + \frac{(L+1)(L-1)(L-2)}{2} \lambda^{L+2}
  + \frac{L(L-1)(L-2)}{6} \lambda^{L+3}}{(1-\lambda^2)^7}.
\nonumber \\
\label{gen_four_vectors_4}
\end{eqnarray}

\section{\label{appD}Molien generating function for five vectors}

The Molien generating function for five vectors can be written as a single rational function:
\begin{equation}
g_\sothree^{\left(a\right)}\left(5,L;\lambda\right) =
\frac{\mathcal{N}_\sothree^{\left(a\right)}\left(5,L;\lambda \right)}{\left(1-\lambda^2\right)^{12}},
\label{gen_five_vectors_1}
\end{equation}
with the numerator equal to:
\begin{eqnarray}
\lefteqn{
\mathcal{N}_\sothree^{\left(a\right)}\left(5, L ; \lambda \right)
=
  \frac{(L+4)(L+3)(L+2)(L+1)}{24} \lambda^L 
+ \frac{(L+4)(L+3)(L+2)L}{6} \lambda^{L+1} }
\nonumber \\
&&
+ \frac{(L+4)(L+3)(L+2)(L+1)}{8} \lambda^{L+2} 
- \frac{(L+4)(L+3)(L^2-3L-\frac{5}{2})}{3} \lambda^{L+3} 
\nonumber \\
&&
- \frac{(L+4)(L+3)(L-3)(7L+2)}{12} \lambda^{L+4} 
- \frac{(L+4)(L-3)(2L+1)}{2} \lambda^{L+5} 
\nonumber \\
&&
+ \frac{(L+4)(L-2)(L-3)(7L+5)}{12} \lambda^{L+6}
+ \frac{(L-2)(L-3)(L^2+5L+\frac{3}{2})}{3} \lambda^{L+7}
\nonumber \\
&&
- \frac{L(L-1)(L-2)(L-3)}{8} \lambda^{L+8}
- \frac{(L+1)(L-1)(L-2)(L-3)}{6} \lambda^{L+9} 
\nonumber \\
&&
- \frac{L(L-1)(L-2)(L-3)}{24} \lambda^{L+10}.
\end{eqnarray}

The six next alternative expressions of the Molien generating function are written as a sum over four rational functions:
\begin{eqnarray}
\lefteqn{g_\sothree^{\left(x\right)}\left(5,L;\lambda \right)}
\nonumber \\
&=&
\frac{\mathcal{N}_\sothree^{\left(x\right),1}\left(5,L;\lambda\right)}{\left(1-\lambda^2\right)^{12}}
+
\frac{\mathcal{N}_\sothree^{\left(x\right),2}\left(5,L;\lambda\right)}{\left(1-\lambda^2\right)^{11}}
+
\frac{\mathcal{N}_\sothree^{\left(x\right),3}\left(5,L;\lambda\right)}{\left(1-\lambda^2\right)^{10}}
+
\frac{\mathcal{N}_\sothree^{\left(x\right),4}\left(5,L;\lambda\right)}{\left(1-\lambda^2\right)^{9}}, \,
x \in \left\{b,c,d,e,f,g\right\},
\nonumber
\end{eqnarray}
and provide for any value of $L$ at least one expression with only non--negative coefficients in the numerators. An heuristic to derive these numerators is presented in section~\ref{S4B2}.

\subsubsection{Numerators of the Molien generating function $(b)$}

\begin{eqnarray}
\lefteqn{\mathcal{N}_\sothree^{\left(b\right),1}\left(5,L;\lambda\right)}
\nonumber \\
&=&
   \frac{(L+4)(L+3)(L+2)(L+1)}{24} \lambda^L 
  + 20(2L+1) \lambda^{L+1}
\nonumber \\
&&  + (-\frac{1}{24}L^4-\frac{5}{12}L^3-\frac{35}{24}L^2+\frac{455}{12}L+19) \lambda^{L+2}
\nonumber \\
\lefteqn{\mathcal{N}_\sothree^{\left(b\right),2}\left(5,L;\lambda \right)}
\nonumber \\
&=&
  (\frac{1}{6}L^4+\frac{3}{2}L^3+\frac{13}{3}L^2-36L-20) \lambda^{L+1} 
\nonumber \\
&&
  - \frac{(L-3)(L^3+13L^2-406L-208)}{24} \lambda^{L+2} 
\nonumber \\
&&
  - \frac{(L-3)(L^3+12L^2-58L-30)}{6} \lambda^{L+3}
\nonumber \\
\lefteqn{\mathcal{N}_\sothree^{\left(b\right),3}\left(5,L;\lambda\right)}
\nonumber \\
&=&
  - \frac{(L-2)(L-3)(L^2-81L-40)}{24} \lambda^{L+2} 
  - \frac{(L-2)(L-3)(L^2-10L-6)}{6} \lambda^{L+3}
\nonumber \\
\lefteqn{\mathcal{N}_\sothree^{\left(b\right),4}\left(5,L;\lambda\right)}
\nonumber \\
&=&
    \frac{L(L-1)(L-2)(L-3)}{4} \lambda^{L+2} 
  + \frac{(L+1)(L-1)(L-2)(L-3)}{6} \lambda^{L+3} 
\nonumber \\
&&
  + \frac{L(L-1)(L-2)(L-3)}{24} \lambda^{L+4}
\nonumber
\end{eqnarray}

\subsubsection{Numerators of the Molien generating function $(c)$}

\begin{eqnarray}
\lefteqn{\mathcal{N}_\sothree^{\left(c\right),1}\left(5,L;\lambda\right)}
\nonumber \\
&=&
   \frac{(L+4)(L+3)(L+2)(L+1)}{24} \lambda^L 
  + 20(2L+1) \lambda^{L+1} 
\nonumber \\
&&
  + (-\frac{1}{24}L^4-\frac{5}{12}L^3-\frac{35}{24}L^2+\frac{455}{12}L+19) \lambda^{L+2}
\nonumber \\
\lefteqn{\mathcal{N}_\sothree^{\left(c\right),2}\left(5,L;\lambda\right)}
\nonumber \\
&=&
    5(2L+1)(2L-7) \lambda^{L+1} 
  - \frac{(L-3)(L^3+13L^2-406L-208)}{24} \lambda^{L+2}
\nonumber \\
\lefteqn{\mathcal{N}_\sothree^{\left(c\right),3}\left(5,L;\lambda\right)}
\nonumber \\
&=&
    \frac{(L-3)(L^3+12L^2-58L-30)}{6} \lambda^{L+1} 
\nonumber \\
&&
  - \frac{(L-2)(L-3)(L^2-81L-40)}{24} \lambda^{L+2} 
\nonumber \\
&&
  - \frac{(L-2)(L-3)(L^2-10L-6)}{6} \lambda^{L+3}
\nonumber \\
\lefteqn{\mathcal{N}_\sothree^{\left(c\right),4}\left(5,L;\lambda\right)}
\nonumber \\
&=&
    \frac{L(L-1)(L-2)(L-3)}{4} \lambda^{L+2} 
  + \frac{(L+1)(L-1)(L-2)(L-3)}{6} \lambda^{L+3} 
\nonumber \\
&&
  + \frac{L(L-1)(L-2)(L-3)}{24} \lambda^{L+4}
\nonumber
\end{eqnarray}

\subsubsection{Numerators of the Molien generating function $(d)$}

\begin{eqnarray}
\lefteqn{\mathcal{N}_\sothree^{\left(d\right),1}\left(5,L;\lambda\right)}
\nonumber \\
&=&
    20(2L+1) \lambda^L 
  + 20(2L+1) \lambda^{L+1}
\nonumber \\
\lefteqn{\mathcal{N}_\sothree^{\left(d\right),2}\left(5,L;\lambda\right)}
\nonumber \\
&=&
    (\frac{1}{24}L^4+\frac{5}{12}L^3+\frac{35}{24}L^2-\frac{455}{12}L-19) \lambda^L 
  + 5(2L+1)(2L-7) \lambda^{L+1} 
\nonumber \\
&&
  - \frac{(L-3)(L^3+13L^2-406L-208)}{24} \lambda^{L+2}
\nonumber \\
\lefteqn{\mathcal{N}_\sothree^{\left(d\right),3}\left(\left( L \right) ;
\Gamma_5 ; \lambda \right)}
\nonumber \\
&=&
    \frac{(L-3)(L^3+12L^2-58L-30)}{6} \lambda^{L+1} 
  - \frac{(L-2)(L-3)(L^2-81L-40)}{24} \lambda^{L+2} 
\nonumber \\
&&
  - \frac{(L-2)(L-3)(L^2-10L-6)}{6} \lambda^{L+3}
\nonumber \\
\lefteqn{\mathcal{N}_\sothree^{\left(d\right),4}\left(5,L;\lambda\right)}
\nonumber \\
&=&
    \frac{L(L-1)(L-2)(L-3)}{4} \lambda^{L+2} 
  + \frac{(L+1)(L-1)(L-2)(L-3)}{6} \lambda^{L+3} 
\nonumber \\
&&
  + \frac{L(L-1)(L-2)(L-3)}{24} \lambda^{L+4}
\nonumber
\end{eqnarray}

\subsubsection{Numerators of the Molien generating function $(e)$}

\begin{eqnarray}
\lefteqn{\mathcal{N}_\sothree^{\left(e\right),1}\left(5,L;\lambda\right)}
\nonumber \\
&=&
    20(2L+1) \lambda^L 
  + 20(2L+1) \lambda^{L+1}
\nonumber \\
\lefteqn{\mathcal{N}_\sothree^{\left(e\right),2}\left(5,L;\lambda\right)}
\nonumber \\
&=&
   (\frac{1}{24}L^4+\frac{5}{12}L^3+\frac{35}{24}L^2-\frac{455}{12}L-19) \lambda^L 
 + 5(2L+1)(2L-7) \lambda^{L+1} 
\nonumber \\
&&
 - \frac{(L-3)(L^3+13L^2-406L-208)}{24} \lambda^{L+2}
\nonumber \\
\lefteqn{\mathcal{N}_\sothree^{\left(e\right),3}\left(5,L;\lambda\right)}
\nonumber \\
&=&
    (L-3)(2L-7)(2L+1) \lambda^{L+1} 
  - \frac{(L-2)(L-3)(L^2-81L-40)}{24} \lambda^{L+2}
\nonumber \\
\lefteqn{\mathcal{N}_\sothree^{\left(e\right),4}\left(5,L;\lambda\right)}
\nonumber \\
&=&
    \frac{(L-2)(L-3)(L^2-10L-6)}{6} \lambda^{L+1} 
  + \frac{L(L-1)(L-2)(L-3)}{4} \lambda^{L+2} 
\nonumber \\
&&
  + \frac{(L+1)(L-1)(L-2)(L-3)}{6} \lambda^{L+3} 
  + \frac{L(L-1)(L-2)(L-3)}{24} \lambda^{L+4}
\nonumber
\end{eqnarray}

\subsubsection{Numerators of the Molien generating function $(f)$}

\begin{eqnarray}
\lefteqn{\mathcal{N}_\sothree^{\left(f\right),1}\left(5,L;\lambda\right)}
\nonumber \\
&=&
    20(2L+1) \lambda^L 
  + 20(2L+1) \lambda^{L+1}
\nonumber \\
\lefteqn{\mathcal{N}_\sothree^{\left(f\right),2}\left(5,L;\lambda\right)}
\nonumber \\
&=&
    5(2L+1)(2L-9) \lambda^L 
  + 5(2L+1)(2L-7) \lambda^{L+1}
\nonumber \\
\lefteqn{\mathcal{N}_\sothree^{\left(f\right),3}\left(5,L;\lambda\right)}
\nonumber \\
&=&
    \frac{(L-3)(L^3+13L^2-406L-208)}{24} \lambda^L 
  + (L-3)(2L+1)(2L-7) \lambda^{L+1} 
\nonumber \\
&&
  - \frac{(L-2)(L-3)(L^2-81L-40)}{24} \lambda^{L+2}
\nonumber \\
\lefteqn{\mathcal{N}_\sothree^{\left(f\right),4}\left(5,L;\lambda\right)}
\nonumber \\
&=&
    \frac{(L-2)(L-3)(L^2-10L-6)}{6} \lambda^{L+1} 
  + \frac{L(L-1)(L-2)(L-3)}{4} \lambda^{L+2} 
\nonumber \\
&&
  + \frac{(L+1)(L-1)(L-2)(L-3)}{6} \lambda^{L+3} 
  + \frac{L(L-1)(L-2)(L-3)}{24} \lambda^{L+4}
\nonumber
\end{eqnarray}

\subsubsection{Numerators of the Molien generating function $(g)$}

\begin{eqnarray}
\lefteqn{\mathcal{N}_\sothree^{\left(g\right),1}\left(5,L;\lambda\right)}
\nonumber \\
&=&
    20(2L+1) \lambda^L
  + 20(2L+1) \lambda^{L+1}
\nonumber \\
\lefteqn{\mathcal{N}_\sothree^{\left(g\right),2}\left(5,L;\lambda\right)}
\nonumber \\
&=&
    5(2L+1)(2L-9) \lambda^L 
  + 5(2L+1)(2L-7) \lambda^{L+1}
\nonumber \\
\lefteqn{\mathcal{N}_\sothree^{\left(g\right),3}\left(5,L;\lambda\right)}
\nonumber \\
&=&
    2(L-3)(L-6)(2L+1) \lambda^L
  + (L-3)(2L+1)(2L-7) \lambda^{L+1}
\nonumber \\
\lefteqn{\mathcal{N}_\sothree^{\left(g\right),4}\left(5,L;\lambda\right)}
\nonumber \\
&=&
    \frac{(L-2)(L-3)(L^2-81L-40)}{24} \lambda^L 
  + \frac{(L-2)(L-3)(L^2-10L-6)}{6} \lambda^{L+1} 
\nonumber \\
&&
  + \frac{L(L-1)(L-2)(L-3)}{4} \lambda^{L+2} 
  + \frac{(L+1)(L-1)(L-2)(L-3)}{6} \lambda^{L+3} 
\nonumber \\
&&
  + \frac{L(L-1)(L-2)(L-3)}{24} \lambda^{L+4}
\nonumber
\end{eqnarray}

\bibliographystyle{spphys}
\bibliography{references}

\begin{thebibliography}{10}
\providecommand{\url}[1]{{#1}}
\providecommand{\urlprefix}{URL }
\expandafter\ifx\csname urlstyle\endcsname\relax
  \providecommand{\doi}[1]{DOI \discretionary{}{}{}#1}\else
  \providecommand{\doi}{DOI \discretionary{}{}{}\begingroup
  \urlstyle{rm}\Url}\fi

\bibitem{wigner}
E.P. Wigner, \emph{Group theory and its application to the quantum mechanics of
  atomic spectra}, \emph{Pure and applied physics}, vol.~5 (Academic Press, New
  York, 1959)

\bibitem{mcweeny}
R.~McWeeny, \emph{Symmetry: an introduction to group theory and its
  applications} (Dover Publications, Mineola, New York, 2002)

\bibitem{hamermesh}
M.~Hamermesh, \emph{Group theory and its application to physical problems}
  (Dover Publications, Mineola, New York, 1989)

\bibitem{dieudonne}
J.A. Dieudonn\'e, J.B. Carrell, Adv. Math. \textbf{4}(1), 1 (1970).
\newblock \doi{10.1016/0001-8708(70)90015-0}

\bibitem{sloane}
N.J.A. Sloane, Amer. Math. Monthly \textbf{84}(2), 82 (1977).
\newblock \doi{10.2307/2319929}

\bibitem{sturmfels}
B.~Sturmfels, \emph{Algorithms in invariant theory}, 2nd edn.
\newblock Texts and monographs in symbolic computation (Springer--Verlag, Wien,
  2008)

\bibitem{doc_00628}
L.~Michel, B.I. Zhilinski\'{\i}, Phys. Rep. \textbf{341}(1--6), 11 (2001)

\bibitem{schmelzer}
A.~Schmelzer, J.N. Murrell, Int. J. Quantum Chem. \textbf{28}(2), 287 (1985).
\newblock \doi{10.1002/qua.560280210}

\bibitem{ischtwan}
J.~Ischtwan, S.D. Peyerimhoff, Int. J. Quantum Chem. \textbf{45}(5), 471
  (1993).
\newblock \doi{10.1002/qua.560450505}

\bibitem{doc_03613}
X.~Huang, B.J. Braams, J.M. Bowman, J. Chem. Phys. \textbf{122}(4), 044308 (12
  pages). (2005).
\newblock \doi{10.1063/1.1834500}.
\newblock {S}ee also the publisher's note, Huang X, Braams B J and Bowman J M
  2007 J. Chem. Phys. 127 099904

\bibitem{doc_03158}
B.J. Braams, J.M. Bowman, Int. Rev. Phys. Chem. \textbf{28}(4), 577 (2009).
\newblock \doi{10.1080/01442350903234923}

\bibitem{doc_03622}
Z.~Xie, J.M. Bowman, J. Chem. Theory Comput. \textbf{6}(1), 26 (2010).
\newblock \doi{10.1021/ct9004917}

\bibitem{a1776}
P.~Cassam-Chena\"{\i}, F.~Patras, J. Math. Chem. \textbf{44}(4), 938 (2008).
\newblock \doi{10.1007/s10910-008-9354-y}

\bibitem{CCDP_methane}
P.~Cassam-Chena\"{\i}, G.~Dhont, F.~Patras, J. Math. Chem. \textbf{53}(1), 58
  (2015).
\newblock \doi{10.1007/s10910-014-0410-5}

\bibitem{kopsky1}
V.~Kopsk\'y, J. Phys. C: Solid State Phys. \textbf{8}(20), 3251 (1975).
\newblock \doi{10.1088/0022-3719/8/20/004}

\bibitem{kopsky2}
V.~Kopsk\'y, J. Phys. A.: Math. Gen. \textbf{12}(4), 429 (1979).
\newblock \doi{10.1088/0305-4470/12/4/004}

\bibitem{rivlin1}
A.C. Pipkin, R.S. Rivlin, Arch. Ration. Mech. An. \textbf{4}(1), 129 (1959).
\newblock \doi{10.1007/BF00281382}

\bibitem{rivlin2}
R.S. Rivlin, Arch. Ration. Mech. An. \textbf{4}(1), 262 (1959).
\newblock \doi{10.1007/BF00281392}

\bibitem{desmorat2020}
B.~Desmorat, M.~Olive, N.~Auffray, R.~Desmorat, B.~Kolev, Computation of
  minimal covariants bases for {2D} coupled constitutive laws (2020).
\newblock \urlprefix\url{https://hal.archives-ouvertes.fr/hal-02888267}

\bibitem{taurines2020}
J.~Taurines, M.~Olive, R.~Desmorat, O.~Hubert, B.~Kolev, Integrity bases for
  cubic nonlinear magnetostriction (2020).
\newblock \urlprefix\url{https://hal.archives-ouvertes.fr/hal-03051787}

\bibitem{sartori}
G.~Sartori, V.~Talamini, J. Math. Phys. \textbf{39}(4), 2367 (1998).
\newblock \doi{10.1063/1.532294}

\bibitem{Planat11}
M.~Planat, Int J. Geom. Methods M. \textbf{08}(02), 303 (2011).
\newblock \doi{10.1142/S0219887811005142}

\bibitem{Holweck14}
F.~Holweck, J.G. Luque, J.Y. Thibon, J. Math. Phys. \textbf{55}(1), 012202
  (2014).
\newblock \doi{10.1063/1.4858336}

\bibitem{Holweck17}
F.~Holweck, J.G. Luque, J.Y. Thibon, J. Math. Phys. \textbf{58}(2), 022201
  (2017).
\newblock \doi{10.1063/1.4975098}

\bibitem{RODRIGUEZBAZAN20211}
E.~{Rodriguez Bazan}, E.~Hubert, J. Symb. Comput. \textbf{107}, 1 (2021).
\newblock \doi{10.1016/j.jsc.2021.01.004}

\bibitem{Li_2013}
J.~Li, B.~Jiang, H.~Guo, J. Chem. Phys. \textbf{139}(20), 204103 (2013).
\newblock \doi{10.1063/1.4832697}

\bibitem{Shao_2016}
K.~Shao, J.~Chen, Z.~Zhao, D.H. Zhang, J. Chem. Phys. \textbf{145}(7), 071101
  (2016).
\newblock \doi{10.1063/1.4961454}

\bibitem{Nandi_2019}
A.~Nandi, C.~Qu, J.M. Bowman, J. Chem. Phys. \textbf{151}(8), 084306 (2019).
\newblock \doi{10.1063/1.5119348}

\bibitem{stanleycov}
R.~Stanley, Proceedings of Symposia in Pure Mathematics \textbf{34}, 345 (1979)

\bibitem{doc_03181}
G.~Dhont, B.I. Zhilinski\'{\i}, J. Phys. A: Math. Theor. \textbf{46}(45),
  455202 (27 pages) (2013)

\bibitem{doc_06124}
G.~Dhont, F.~Patras, B.I. Zhilinski\'{\i}, J. Phys. A: Math. Theor.
  \textbf{48}(3), 035201 (19 pages) (2015).
\newblock \doi{10.1088/1751-8113/48/3/035201}

\bibitem{doc_08519}
G.~Dhont, B.I. Zhilinski{\'i}, in \emph{Geometric Methods in Physics: XXXIV
  Workshop, Bia{\l}owie{\.{z}}a, Poland, June 28 -- July 4, 2015}, ed. by
  P.~Kielanowski, T.S. Ali, P.~Bieliavsky, A.~Odzijewicz, M.~Schlichenmaier,
  T.~Voronov (Springer International Publishing, Cham, 2016), pp. 105--114.
\newblock \doi{10.1007/978-3-319-31756-4\_11}

\bibitem{doc_01743}
T.~Molien, Sitzungsber. K{\"o}nig. Preuss. Akad. Wiss. \textbf{52}, 1152 (1897)

\bibitem{sutcliffe}
B.T. Sutcliffe, J. Chem. Soc., Faraday Trans. \textbf{89}(14), 2321 (1993).
\newblock \doi{10.1039/FT9938902321}

\bibitem{GATTI20091}
F.~Gatti, C.~Iung, Phys Rep. \textbf{484}(1), 1 (2009).
\newblock \doi{10.1016/j.physrep.2009.05.003}

\bibitem{VMK}
D.A. Varshalovich, A.N. Moskalev, V.K. Khersonskii, \emph{Quantum theory of
  angular momentum} (World Scientific, Singapore, 1988)

\bibitem{a1610}
M.A. Collins, D.F. Parsons, J. Chem. Phys. \textbf{99}(9), 6756 (1993)

\bibitem{b100}
H.~Weyl, \emph{The classical groups. Their invariants and representations}
  (Princeton University Press, Princeton, New Jersey, 1939)

\bibitem{HR1974}
M.~Hochster, J.L. Roberts, Adv. Math. \textbf{13}(2), 115 (1974).
\newblock \doi{10.1016/0001-8708(74)90067-X}

\bibitem{DK}
H.~Derksen, G.~Kemper, \emph{Computational Invariant Theory} (Springer Berlin
  Heidelberg, Berlin, Heidelberg, 2015).
\newblock \doi{10.1007/978-3-662-48422-7}

\bibitem{VanDenBergh}
M.~{Van den Bergh}, in \emph{Proceedings of the International Congress of
  Mathematicians}, ed. by S.D. Chatterji (Birkh{\"a}user Basel, Basel, 1995),
  pp. 352--362.
\newblock \doi{10.1007/978-3-0348-9078-6\_29}

\bibitem{Stanley82}
R.P. Stanley, Invent. Math. \textbf{68}(2), 175 (1982)

\bibitem{Schmiedt15}
H.~Schmiedt, S.~Schlemmer, P.~Jensen, J. Chem. Phys. \textbf{143}(15), 154302
  (2015).
\newblock \doi{10.1063/1.4933001}

\bibitem{herzberg}
G.~Herzberg, Nature \textbf{163}, 170 (1949).
\newblock \doi{10.1038/163170a0}

\bibitem{camparguePRR}
A.~Campargue, S.~Kassi, A.~Yachmenev, A.A. Kyuberis, J.~K\"upper, S.N.
  Yurchenko, Phys. Rev. Research \textbf{2}, 023091 (2020).
\newblock \doi{10.1103/PhysRevResearch.2.023091}

\bibitem{camparguePCCP}
A.~Campargue, A.M. Solodov, A.A. Solodov, A.~Yachmenev, S.N. Yurchenko, Phys.
  Chem. Chem. Phys. \textbf{22}, 12476 (2020).
\newblock \doi{10.1039/D0CP01667E}

\bibitem{yachmenev}
A.~Yachmenev, A.~Campargue, S.N. Yurchenko, J.~K\"upper, J.~Tennyson, J. Chem.
  Phys. \textbf{154}(21), 211104 (2021).
\newblock \doi{10.1063/5.0053279}

\bibitem{BL}
L.C. Biedenharn, J.D. Louck, \emph{Angular momentum in quantum physics},
  \emph{Encyclopedia of mathematics and its applications}, vol.~8 (Cambridge
  University Press, 1985)

\bibitem{equation_3_616_7_gradshteyn_ryzhik}
I.S. Gradshteyn, I.M. Ryzhik, \emph{Table of integrals, series, and products}
  (Academic Press, New York, 1980).
\newblock Eq.~(3.616.7)

\bibitem{BL_3.153}
{Eq.}~(3.153) of Ref.~\cite{BL}

\bibitem{BL_6.183}
{Eq.}~(6.183) of Ref.~\cite{BL}

\bibitem{steinborn}
E.O. Steinborn,  (Academic Press, 1973), pp. 83--112.
\newblock \doi{10.1016/S0065-3276(08)60559-6}.
\newblock {Eq.}~(123)

\bibitem{BFB}
M.A. Blanco, M.~Fl\'orez, M.~Bermejo, J. Mol. Struct. THEOCHEM \textbf{419}(1),
  19 (1997).
\newblock \doi{10.1016/S0166-1280(97)00185-1}.
\newblock {Eq.}~(6)

\end{thebibliography}

\end{document}